\documentclass[openacc]{rsproca_new}
\usepackage{psfrag}
\usepackage{subfigure}

\def\H{{\mathcal{H}}}
\def\O{{\mathcal{O}}}

\def\P{{\mathcal{P}}}
\def\Q{{\mathcal{Q}}}
\def\R{{\mathcal{R}}}
\def\E{{\mathcal{E}}}

\def\K{{\mathcal{K}}}
\def\L{{\mathcal{L}}}

\def\X{{\mathcal{X}}}
\def\W{{\mathcal{W}}}

\def\RR{{\mathbb{R}}}

\def\TT{{\mathbb{T}}}
\def\ZZ{{\mathbb{Z}}}

\def\e{{\mathbf e}}

\def\x{{\bm x}}

\def\u{{\bm u}}

\def\x{{\mathbf x}}

\def\z{{\mathbf z}}
\def\f{{\mathbf f}}

\def\u{{\mathbf u}}
\def\v{{\mathbf v}}

\def\0{{\mathbf 0}}

\def\bomega{\boldsymbol{\omega}}
\def\bnabla{\boldsymbol{\nabla}}
\def\bDelta{\boldsymbol{\Delta}}

\def\Dpartial#1#2{ {\partial #1 \over \partial #2} }

\def\Bmp#1{ \begin{minipage}{#1} }
\def\Emp{ \end{minipage} }
\def\Bmpc#1{ \begin{minipage}[c]{#1} }
\def\Bmpt#1{ \begin{minipage}[t]{#1} }
\def\Bmpb#1{ \begin{minipage}[b]{#1} }

\def\tTE{\widetilde{T}_{\E_0}}

\newcommand{\xvec}{\mathbf{x}}
\newcommand{\uvec}{\mathbf{u}}

\newcommand{\laplacian}{\Delta}
\newcommand{\rot}{\bnabla\times}

\newcommand{\tomega}{\widetilde{\omega}}
\newcommand{\tuE}{\widetilde{\mathbf{u}}_{\E_0}}

\newcommand{\tuET}{\widetilde{\uvec}_{0;\E_0,T}}

\newcommand{\tuBT}{\widetilde{\uvec}_{0;B,T}}

\newcommand{\tuEtT}{\widetilde{\uvec}_{0;\E_0,\tTE}}

\newcommand{\ET}{\E_T(\u_0)}

\newcommand{\tuuE}{\widetilde{u}_{\E_0}}
\newcommand{\tuuET}{\widetilde{u}_{0;\E_0,T}}
\newcommand{\tpsiKP}{\widetilde{\psi}_{\K_0,\P_0}}

\newcommand{\argmax}{\operatorname{argmax}}

\newcommand{\Id}{\operatorname{Id}}

\newtheorem{problem}{Problem}

\begin{document}

\title{Systematic Search For Extreme and Singular Behavior in Some Fundamental Models of Fluid Mechanics}

\author{
B. Protas$^{1}$}

\address{$^{1}$Department of Mathematics and Statistics, McMaster University,
  Hamilton, ON, Canada}

\subject{Fluid Mechanics, Applied Mathematics}


\corres{Bartosz Protas\\
\email{bprotas@mcmaster.ca}}

\begin{abstract}
  This review article offers a survey of the research program focused
  on a systematic computational search for extreme and potentially
  singular behavior in hydrodynamic models motivated by open questions
  concerning the possibility of a finite-time blow-up in the solutions
  of the Navier-Stokes system. Inspired by the seminal work of Lu \&
  Doering (2008), we sought such extreme behavior by solving PDE
  optimization problems with objective functionals chosen based on
  certain conditional regularity results and a priori estimates
  available for different models. No evidence for singularity
  formation was found in extreme Navier-Stokes flows constructed in
  this manner in 3D.  We also discuss the results obtained for 1D
  Burgers and 2D Navier-Stokes systems, and while singularities are
  ruled out in these flows, the results presented provide interesting
  insights about sharpness of different energy-type estimates known
  for these systems. Connections to other bounding techniques are also
  briefly discussed.
\end{abstract}


\begin{fmtext}
\section{Introduction}
\label{sec:intro}

One of the central problems in mathematical fluid mechanics is the
question whether the Navier-Stokes system, which is the main
mathematical model used to describe the motion of viscous
incompressible fluids, admits unique classical solutions existing
globally in time for all sufficiently regular initial data
\cite{d09,Robinson2020}. In other words, the question is whether
starting from such smooth initial data it may be possible for a
``singularity'' to form spontaneously in the solution such that the
equations would no longer be satisfied in the classical pointwise
sense. Should such situation indeed occur, this would 
\end{fmtext}

\maketitle

\noindent
invalidate the Navier-Stokes system as an acceptable model to describe
flows of viscous incompressible fluids. While for Navier-Stokes flows
in two dimensions (2D) the problem is solved and finite-time blow-up
has been ruled out \cite{kl04}, in the physically more relevant case
of three-dimensional (3D) flows the problem remains open. Recognizing
the difficulty and significance of this problem, the Clay Mathematics
Institute named it one of its seven ``millennium problems'' posed as
challenges to the mathematical community at the beginning of the 21st
century \cite{f00}.  On the other hand, weak solutions, which may in
principle be nonunique and involve singularities, are known to exist
globally in time since the work of Leray \cite{l34} and nonuniqueness
was recently established in \cite{BuckmasterVicol2019} for weak
solutions of a certain type.  Analogous questions concerning existence
of unique smooth solutions also remain open for the inviscid Euler
equation in 3D \cite{gbk08}.

When fundamental properties of its solutions are studied, the
Navier-Stokes system is usually considered on domains without solid
boundaries, namely, the unbounded domain $\Omega = \RR^3$ or a
periodic box (3D torus) $\Omega = {\TT_L^3 := [0,L]^3}$,
{where $L>0$ is the domain size,} and here we will focus on the
latter case.  Assuming we are interested in solutions on the time
interval $[0,T]$, the Navier-Stokes system is defined as
\begin{subequations}\label{eq:NS}
\begin{alignat}{2}
\partial_t\u + \u\cdot\bnabla\u + \bnabla p - \nu\laplacian\u & = 0 & &\qquad\mbox{in} \,\,\Omega\times(0,T], \label{eq:NSa}  \\
\bnabla\cdot\u & = 0 & & \qquad\mbox{in} \,\,\Omega\times[0,T], \label{eq:NSb} \\
\u(0) & = \u_0, &   &
\end{alignat}
\end{subequations}
where $\u \; :\ \; [0,T] \times \Omega \rightarrow \RR^3$ and $p \; :\
\; [0,T] \times \Omega \rightarrow \RR$ are the velocity and pressure
fields, $\nu>0$ is the coefficient of kinematic viscosity, $\u_0$ is a
divergence-free initial condition whereas the fluid density $\rho$ is
assumed constant and equal to unity ($\rho\equiv 1$).  Without loss of
generality, we will assume the initial data $\u_0$ to have zero mean.

{In system \eqref{eq:NS} there are three physical parameters:
  the domain size $L$, kinetic viscosity $\nu$ and the ``magnitude''
  of the initial data $\u_0$. They can be combined into a single
  dimensionless quantity, the Reynolds number, meaning that only one of
  these parameters needs to be changed in order to study solutions of
  \eqref{eq:NS} in different regimes. In the investigations surveyed
  here one typically considers variations of the size of the initial
  data $\u_0$ while fixing $L$ and $\nu$. Thus, in keeping with these
  earlier studies, we will henceforth set $L = 1$ and define $\TT :=
  \TT_1$, however, explicit dependence on $\nu$ will be retained in
  some of the estimates.}

Important quantities characterizing solutions of system \eqref{eq:NS}
include the Lebesgue norms of the velocity field
\begin{equation}
\|\u(t))\|_{L^q(\Omega)}  := \left( \int_\Omega |\u(t,\x)|^q \,d\x \right)^{\frac{1}{q}}, \qquad q \ge 1,
\label{eq:uLq}
\end{equation}
as well as the kinetic energy and enstrophy\footnote{{We note that
    unlike energy, cf.~\eqref{eq:K}, enstrophy is often defined
    without the factor of 1/2. However, for consistency with earlier
    studies belonging to this research program
    \cite{ap11a,ap13a,ap13b,ap16,Yun2018,KangYunProtas2020}, we choose
    to retain this factor here.}} defined as
\begin{subequations}
\begin{align}
\K(\u(t))  & := \frac{1}{2} \|\u(t))\|_{L^2(\Omega)}, \label{eq:K} \\
\E(\u(t))  & :=  \frac{1}{2}\int_\Omega | \bomega(t,\x) |^2 \,d\x = \frac{1}{2} \|\bnabla\u(t))\|_{L^2(\Omega)}, \label{eq:E}
\end{align}
\end{subequations}
where $\bomega(t,\x) := \rot\u(t,\x)$ is the vorticity (``$:=$'' means
``equal to by definition''). {In addition, we will also use
  Sobolev spaces} $H^s(\Omega)$, $s \in \RR^+$, of functions with
square-integrable weak derivatives of order $s$ \cite{af05}.

The question about the possibility of singularity formation in
solutions of the Navier-Stokes system \eqref{eq:NS} is primarily a
problems in mathematical analysis of partial differential equations
(PDEs).  An important class of results obtained to date has the form
of ``conditional regularity results'' stating conditions which need to
be satisfied by a Leray-Hopf weak solution for it to also satisfy
system \eqref{eq:NS} in the classical sense, i.e., pointwise in $(0,T]
\times \Omega$. Typically, such solutions will then also be smooth
(real-analytic) \cite{RobinsonRodrigoSadowski2016}. Conditional
regularity results are often accompanied by a priori estimates
involving some related quantities and also applicable to weak
solutions.  Arguably, the best known conditional regularity result is
the enstrophy condition \cite{ft89} asserting that $\u(t)$ is a smooth
solution of system \eqref{eq:NS} on the time interval $[0,T]$ provided
its enstrophy \eqref{eq:E} remains bounded, i.e.,
\begin{equation}
\label{eq:supE}
\mathop{\sup}_{0 \leq t \leq T} \E(\u(t))  < \infty.
\end{equation}
While it is not known whether \eqref{eq:supE} is true for all initial
data $\u_0$ and arbitrarily large $T$, Leray-Hopf weak solutions
satisfy $\int_0^T \E(\u(t)) \, dt < \infty$ (however, the boundedness
of $\int_0^T \E(\u(t))^2 \, dt$ is an open question).

Another important conditional regularity result is given by the family
of the Ladyzhenskaya-Prodi-Serrin conditions asserting that
Navier-Stokes flows $\u(t)$ are smooth and satisfy system
\eqref{eq:NS} in the classical sense provided that
\cite{KisLad57,Prodi1959,Serrin1962}
\begin{equation}
\u \in L^p([0,T];L^q(\Omega)), \quad 2/p+3/q = 1, \quad q > 3.
\label{eq:LPS}
\end{equation}
These conditions were recently generalized in \cite{Gibbon2018} to
include norms of the derivatives of the velocity field. As regards the
limiting case with $q = 3$, the corresponding condition was
established in \cite{Escauriaza2003}
\begin{equation}
\u \in L^{\infty}([0,T];L^3(\Omega))
\label{eq:LPS3}
\end{equation}
and a related blow-up criterion was recently obtained in
\cite{Tao2020}.  Condition \eqref{eq:LPS} implies that should a
singularity form in a classical solution $\u(t)$ of the Navier-Stokes
system \eqref{eq:NS} at some finite time $0 < t_0 < \infty$, then
necessarily
\begin{equation}
\lim_{t \rightarrow t_0} \int_0^t \| \u(\tau) \|_{L^q(\Omega)}^p \, d\tau \rightarrow \infty, \quad 2/p+3/q = 1, \quad q > 3.
\label{eq:LPSblowup}
\end{equation}
At the same time, the time evolution of the solution norm $\| \u(t)
\|_{L^q(\Omega)}$ on the time interval $[0,T]$ is subject to the some
a priori bounds valid also for Leray-Hopf weak solutions
\cite{Gibbon2018}, which might involve singularities.  An estimate of
this type was discussed in \cite{Constantin1991} and was rederived
with an upper bound explicitly depending on the initial data in
\cite{KangProtas2021}
\begin{equation}
\int_0^T \| \u(\tau) \|_{L^q(\Omega)}^{\frac{4q}{3(q-2)}} \, d\tau \le  C \, \K_0^{\frac{2q}{3(q-2)}}, \qquad 2 \le q \le 6,
\label{eq:LPSbound}
\end{equation}
where $\K_0 := \K(\u_0)$ and $C>0$ is a generic constant whose
numerical value may vary between different estimates. We note that
the integrals in \eqref{eq:LPSblowup} and \eqref{eq:LPSbound} differ
in the exponent in the integrand expressions which is smaller in the
latter case.

{We add that in the context of the inviscid Euler system a
  conditional regularity result analogous to \eqref{eq:supE} and
  \eqref{eq:LPS}--\eqref{eq:LPS3} is given by the Beale-Kato-Majda
  (BKM) criterion which asserts that an Euler flow remains smooth on
  $[0,T]$ if and only if $\int_0^T \| \bomega(\tau)
  \|_{L^\infty(\Omega)} \, d\tau < \infty$ \cite{bkm84}.  A relation
  between potential blow-up in Euler flows and Navier-Stokes flows
  with sufficiently small viscosity was established in
  \cite{Constantin1986}.  Recently, finite-time singularity formation
  in 3D axisymmetric Euler flows on domains exterior to a boundary
  with conical shape was proved in \cite{ElgindiJeong2018}.}

In order to obtain insights about the enstrophy condition
\eqref{eq:supE}, we assume the Navier-Stokes system \eqref{eq:NS}
admits a smooth classical solution $\u(t)$ for times $t \in [0,T]$,
where $T$ is sufficiently small, which is guaranteed by local
existence theorems \cite{RobinsonRodrigoSadowski2016}. We then
consider the equations for the evolution of the kinetic energy
\eqref{eq:K} and the enstrophy \eqref{eq:E} obtained multiplying
\eqref{eq:NSa} by, respectively, $\u$ and $\bDelta\u$, integrating
over $\Omega$ and performing integrations by parts (these operations
are justified for $t \in [0,T]$ since the solution $\u(t)$ is smooth
there)
\begin{subequations}
\begin{align}
\frac{d\K(\u(t))}{dt} & = - \nu \E(\u(t)), \label{eq:dKdt} \\
\frac{d\E(\u(t))}{dt} & = -\nu\int_\Omega|\bDelta\u|^2\,d\x + \int_{\Omega} \u\cdot\bnabla\u\cdot\bDelta\u\, d\x =: \R_\E(\u).
\label{eq:dEdt} 
\end{align}
\end{subequations}
As shown in \cite{ld08}, relation \eqref{eq:dEdt} ca be used to obtain
the following upper bound on the rate of growth of enstrophy
\begin{equation}
\frac{d\E}{dt} \leq \frac{27}{8\,\pi^4\,\nu^3} \E^3. 
\label{eq:dEdt_estimate_E}
\end{equation} 
By simply integrating the differential inequality in
\eqref{eq:dEdt_estimate_E} with respect to time we obtain the
finite-time bound
\begin{equation}
\E(\u(t)) \leq \frac{\E_0}{\sqrt{1 - \frac{27}{4\,\pi^4\,\nu^3}\,\E_0^2\, t}}, 
\label{eq:Et_estimate_E0}
\end{equation}
where $\E_0 := \E(\u_0)$, which becomes infinite at time $t_0 =
4\,\pi^4\,\nu^3 / (27\, \E_0^2)$. Thus, based on inequality
\eqref{eq:Et_estimate_E0}, which is the best estimate available to
date, it is not possible to establish the boundedness of the enstrophy
$\E(\u(t))$ required in condition \eqref{eq:supE} and hence also the
regularity of solutions globally in time. However, boundedness of
enstrophy and existence of smooth solutions can be established for
arbitrarily long times provided the initial data $\u_0$ is ``small'',
more precisely, when $\K_0 \E_0 = \O(\nu^4)$ \cite{ld08}.

In a similar vein, the Ladyzhenskaya-Prodi-Serrin condition
\eqref{eq:LPSblowup} can be studied by considering the rate of growth
of the $L^q$ norm of the velocity field, for which an upper bound was
already known to Leray \cite{l34}, see also
\cite{Giga1986,RobinsonSadowskiSilva2012,RobinsonSadowski2014},
\begin{equation}
\frac{1}{q} \frac{d}{dt} \| \u(t) \|_{L^q(\Omega)}^q \le C  \| \u(t) \|_{L^q(\Omega)}^{\frac{q(q-1)}{q-3}}, \qquad q > 3.
\label{eq:dLqdt}
\end{equation}
However, as was the case with the enstrophy condition, this approach
does not lead to estimates that would allow one to ascertain the
finiteness of the integral expression in \eqref{eq:LPSblowup}.


{While the blow-up problem is fundamentally a question in
  mathematical analysis, a lot of computational studies have been
  carried out since the mid-'80s in order to shed light on the
  hydrodynamic mechanisms which might lead to singularity formation in
  finite time. Given that such flows evolving near the edge of
  regularity involve formation of very small flow structures, these
  computations typically require the use of state-of-the-art
  computational resources available at a given time. The computational
  studies focused on the possibility of finite-time blow-up in the 3D
  Navier-Stokes and/or Euler system include
  \cite{bmonmu83,ps90,b91,k93,p01,bk08,oc08,o08,ghdg08,
    gbk08,h09,opc12,bb12,opmc14,CampolinaMailybaev2018}, all of which
  considered problems defined on domains periodic in all three
  dimensions. The investigations \cite{dggkpv13,k13,gdgkpv14,k13b}
  focused on the time evolution of vorticity moments and compared it
  against bounds on these quantities obtained using rigorous analysis.
  Recent computations \cite{Kerr2018} considered a ``trefoil''
  configuration meant to be defined on an unbounded domain (although
  the computational domain was always truncated to a finite periodic
  box). A simplified semi-analytic model of vortex reconnection was
  recently developed and analyzed based on the Biot-Savart law and
  asymptotic techniques \cite{MoffattKimura2019a,MoffattKimura2019b}.
  We also mention the studies \cite{mbf08} and \cite{sc09}, along with
  references found therein, in which various complexified forms of the
  Euler equation were investigated. The idea of this approach is that,
  since the solutions to complexified equations have singularities in
  the complex plane, singularity formation in the real-valued problem
  is manifested by the collapse of the complex-plane singularities
  onto the real axis.  Overall, the outcome of these investigations is
  rather inconclusive: while for the Navier-Stokes {system most of
    the} recent computations do not offer support for finite-time
  blow-up, the evidence appears split in the case of the Euler system.
  In particular, the studies \cite{bb12} and \cite{opc12} hinted at
  the possibility of singularity formation in finite time. In this
  connection we also highlight the {computational} investigations
  \cite{lh14a,lh14b} in which blow-up was {documented} in axisymmetric
  Euler flows on a bounded (tubular) domain. Recently, numerical
  evidence for blow-up in solutions of the Navier-Stokes system in 3D
  axisymmetric geometry with a degenerate variable diffusion
  coefficient was provided in \cite{HouHuang2021}.}

An entirely different approach designed to {\em systematically} search
for potentially singular Navier-Stokes flows was proposed by Lu \&
Doering based on the conditional regularity result \eqref{eq:supE} in
\cite{ld08} and was later developed in
\cite{ap11a,ap13a,ap16,Yun2018,KangYunProtas2020,KangProtas2021}.  The
idea is to look for initial data which might potentially lead to a
finite-time singularity as a solution of a certain variational
optimization problem with the objective functional and constraints
motivated by estimates
\eqref{eq:dEdt_estimate_E}--\eqref{eq:Et_estimate_E0}. In addition, in
this framework it is also possible to check (usually at the level of
computational evidence) the sharpness of a priori estimates such as
\eqref{eq:dEdt_estimate_E}. We say that a polynomial upper bound of
the type $C \E^{\alpha}$ for some $\alpha > 0$ is sharp (up to a
numerical prefactor) if the expression on the left-hand side (LHS) of
the estimate is $\O(\E^{\alpha})$ as $\E \rightarrow \infty$. A family
of initial conditions and the corresponding flows parameterized by
$\E_0$ and saturating a certain estimate in the above sense is
referred to as ``extreme''. Energy-type estimates similar to
\eqref{eq:dEdt_estimate_E}--\eqref{eq:Et_estimate_E0} are also known
for the one-dimensional (1D) viscous Burgers equation and the
two-dimensional (2D) Navier-Stokes system. While these two systems are
known to be globally well-posed \cite{kl04}, the question whether
these estimates are sharp is in fact quite pertinent, because they are
established using similar mathematical techniques as
\eqref{eq:dEdt_estimate_E}--\eqref{eq:Et_estimate_E0}.  These
estimates are obtained from the governing equations applying different
functional inequalities and although each of these inequalities is
known to be sharp, sharpness need not be preserved if they are chained
together (because different inequalities are saturated by different
fields).

These observations have motivated a research program focused on
probing the sharpness of a number of key estimates, both instantaneous
as in \eqref{eq:dEdt_estimate_E} and finite-time as in
\eqref{eq:Et_estimate_E0}, in the 1D Burgers and 2D Navier-Stokes
systems, in addition to examining estimates
\eqref{eq:dEdt_estimate_E}--\eqref{eq:Et_estimate_E0} and more
recently the Ladyzhenskaya-Prodi-Serrin criterion \eqref{eq:LPSblowup}
in 3D Navier-Stokes flows.  Progress in this research program was
largely enabled by the development of robust computational approaches
for the solution of large-scale PDE-constrained optimization problems.
Since a number of important milestones has recently been attained in
this research program, the present review paper aims to survey these
developments. 

{Most of the optimization problems considered here are
  nonconvex, hence their solutions found numerically based on local
  optimality conditions are local maximizers only. Thus, unless stated
  otherwise, when we refer to ``maximizing solutions'' defined with
  $\mathop{\arg\max}$ we will in fact mean local maximizers.
  Theoretical results concerning existence of (possibly nonunique)
  solutions to optimization problem involving different hydrodynamic
  PDE models are available in the literature, which includes the
  seminal study \cite{control:abergel1} and the monographs
  \cite{Fursikov2000,g03,Troltzsch2010}.}

The structure of the paper is as follows: in the next three sections
we review instantaneous and finite-time energy-type estimates known
for the 1D Burgers, 2D and 3D Navier-Stokes systems, and discuss
different optimization problems that have been introduced to test
their sharpness before presenting some key results (in Section
\ref{sec:1D} devoted to the 1D Burgers equation we also provide
comments about the corresponding stochastic problem and the system
with fractional dissipation); in Section \ref{sec:bound} we draw some
connections to other research problems concerned with establishing
bounds on the behavior of hydrodynamic models such as the background
method and the methods based on sum-of-squares (SOS) polynomial
bounds; finally, summary and conclusions are deferred to Section
\ref{sec:final} where we also provide an outlook; more technical
material concerning the numerical solution of the optimization
problems studied in this research program is collected in an appendix.

\section{Estimates for the 1D Burgers Equation}
\label{sec:1D}

The 1D viscous Burgers equation
\begin{subequations}
\label{eq:Burgers}
\begin{align}
\partial_{t}u+\frac{1}{2}\,\partial_{x}u^{2}-\nu\partial_{xx} u=0 \quad & \mbox{in}\ (0,T]\times \Omega,\label{eq:Burgersa}  \\
u(0) = u_0&,\label{eq:Burgersc}
\end{align}
\end{subequations}
where $\Omega = \TT$ is a periodic domain and $u_0$ the initial
condition, has often been used as a highly idealized model of the
Navier-Stokes system \cite{bk07}. Unlike its inviscid variant
(obtained by setting $\nu = 0$ in \eqref{eq:Burgersa}) which exhibits
a well-documented finite-time blow-up, system \eqref{eq:Burgers} is
globally well posed in the classical sense \cite{kl04}. As shown in
\cite{ld08}, defining the 1D equivalent of enstrophy as
\begin{equation}
\E(u(t))  :=  \frac{1}{2}\int_0^1 \left| \partial_x u(t,x)  \right|^2 \,dx, 
\label{eq:E1}
\end{equation}
it is possible to obtain an estimate for the rate of growth of
enstrophy analogous to \eqref{eq:dEdt_estimate_E} in the
form\footnote{{Due to the presence of the factor 1/2 in
    \eqref{eq:E1}, the coefficients in \eqref{eq:dEdtbound1D} and
    \eqref{eq:r1d} differ from those given in \cite{ld08}. For the
    same reason, relations (4) and (5) in \cite{ap11a} contain
    incorrect prefactors. The second term on the RHS.in relation
    \eqref{eq:r} appears with an incorrect sign in \cite{ld08} and in
    \cite{ap11a}. }}
\begin{equation}
\frac{d\E}{dt} \leq { 3 \left(\frac{1}{2 \pi^2\nu}\right)^{1/3}}\E^{5/3}.
\label{eq:dEdtbound1D}
\end{equation}
Based on this estimate, the corresponding finite-time bound on
enstrophy was obtained in \cite{ap11a} by integrating
\eqref{eq:dEdtbound1D} in time
\begin{equation}
\max_{t \in [0,T]} \E(u(t)) \leq \left[\E_0^{1/3} + {\frac{1}{4}\left(\frac{1}{2\pi^2 \nu}\right)^{4/3}}\E_0\right]^{3} \ \underset {\E_0 \rightarrow \infty} {\longrightarrow} \
{\frac{1}{64} \left(\frac{1}{2 \pi^2 \nu}\right)^{4}} \E_0^{3}.
\label{eq:Etbound1D}
\end{equation}
We emphasize that in contrast to \eqref{eq:Et_estimate_E0} this bound
is valid uniformly in $t$. Moreover, it also exhibits a well-defined
asymptotic behavior in the large enstrophy limit. In this context one
should also mention Biryuk's work \cite{Biryuk2001} which implies a
finite-time bound with a smaller exponent, namely, $\max_{t \ge 0}
\E(u(t)) \le C_B \E_0^{3/2}$. While this approach did not rely on time
integration of an instantaneous estimate such as
\eqref{eq:dEdtbound1D}, the prefactor in this estimate $C_B = C_B(\|
u_0 \|_{H^2})$ requires the $H^2$ norm of the initial data $u_0$ to be
bounded. {Consequently, owing to Poincar\'e's inequality, this
  prefactor will not remain bounded in the limit we are interested in,
  i.e., as $\E_0 \rightarrow \infty$.}

Lu \& Doering \cite{ld08} posed an interesting question about the
sharpness of estimate \eqref{eq:dEdtbound1D} and to elucidate it
formulated the following optimization problem
\begin{problem}
\label{pb:maxdEdt1D}
  Given $\E_0\in\RR_+$ and the objective functional, cf.~\eqref{eq:dEdt},
\begin{equation}
{r(u) := -\nu \| \partial_{xx} u \|_{L^2([0,1])}^2 
- \frac{1}{2} \int_0^1 \left(\partial_x u \right)^3 \, dx = \frac{d\E(u(t))}{dt},}
\label{eq:r}
\end{equation}
find
\begin{equation*}
\tuuE  =  \mathop{\arg\max}_{u \in \Sigma_{\E_0}} \, r(u), \qquad \text{where} \qquad
\Sigma_{\E_0}  :=  \left\{u\in H^2(\Omega)\,\colon \; \int_0^1 u \, dx = 0, \ E(u) = \E_0 \right\}.
\end{equation*} 
\end{problem}
\noindent 
The idea behind this problem is to maximize the LHS in estimate
\eqref{eq:dEdtbound1D} for a range of values of the constraint $\E_0$
to see whether the maximum attainable values of $d\E/dt$,
cf.~\eqref{eq:r}, saturate the upper bound on the RHS, in the sense of
having the same dependence on $\E_0$. Remarkably, Lu \& Doering were
able to solve Problem \ref{pb:maxdEdt1D} in a closed form using the
method of Lagrange multipliers with the optimal solution $\tuuE$
expressed in terms of elliptic integrals and Jacobi elliptic
functions. By analyzing the asymptotic behavior of these solutions for
large enstrophies, they concluded that
\begin{equation}
r(\tuuE) \sim \frac{{0.393}}{\nu^{1/3}} \E_0^{5/3} \quad \text{as} \quad \E_0 \rightarrow \infty,
\label{eq:r1d}
\end{equation}
thus demonstrating that estimate \eqref{eq:dEdtbound1D} is sharp (up
to a numerical prefactor which is larger than in \eqref{eq:r1d} by
about 2.83). In other words, for each value of $\E_0$, the optimal
fields $\tuuE$, which have the form of steep waves with fronts
becoming sharper as $\E_0$ increases, instantaneously produce as much
enstrophy $r(\tuuE)$ as is only allowed by the mathematically rigorous
analysis of the 1D Burgers system \eqref{eq:Burgers}. {On the
  other hand, solving the Burgers system with optimizers $\tuuE$ of
  Problem \ref{pb:maxdEdt1D} used as the initial data produces maximum
  enstrophy which scales as $\O(\E_0)$ for large $\E_0$, far below
  what is allowed by estimate \eqref{eq:Etbound1D}.}

\begin{figure}[t]
\begin{center}
\mbox{\hspace*{-0.35cm}
\subfigure[]{\includegraphics[width=0.5\textwidth]{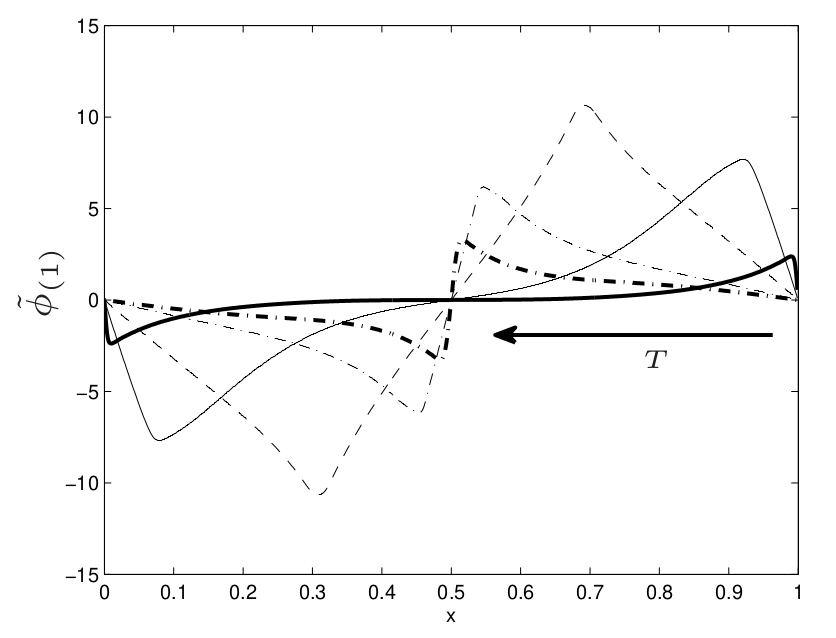}}\quad
\subfigure[]{\includegraphics[width=0.5\textwidth]{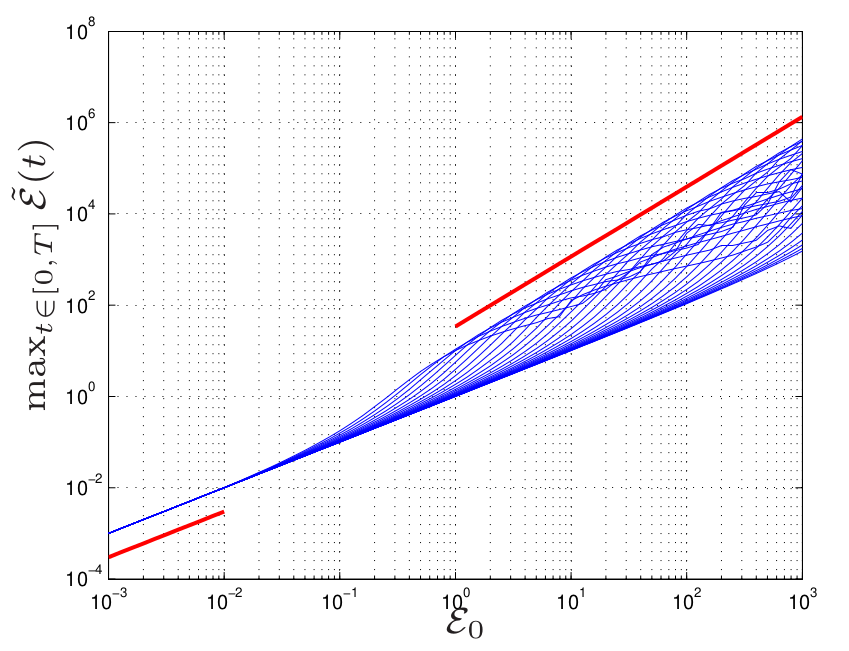}}}
\caption{(a) Optimal initial conditions $\tuuET$ obtained by solving
  Problem \ref{pb:maxET1D} with fixed enstrophy $\E_0 = 10^3$ and
  different time intervals: (thick solid line) $T=10^{-3}$, (thin
  solid line) $T=10^{-2}$, (thin dashed line) $T=10^{-1.5}$, (thin
  dotted line) $T=10^{-1}$ and (thick dotted line) $T=10^{0}$; the
  arrow indicates the trend with increasing $T$.  (b) Maximum
  enstrophy $\max_{t\in[0,T]} \E_T(\tuuET)$ as a function of initial
  enstrophy $\E_0$ for different $T$. Two distinct power laws can be
  observed with exponents 1 for small $\E_0$ and 3/2 for large $\E_0$,
  cf.~\eqref{eq:maxETap}.}
\label{fig:maxET1D}
\vspace*{-1.0cm}
\end{center}
\end{figure}

The companion question about sharpness of the corresponding
finite-time estimate \eqref{eq:Etbound1D} was taken up by Ayala \&
Protas in \cite{ap11a} where the following optimization problem was
considered
\begin{problem}
\label{pb:maxET1D}
  Given $\E_0, T\in\RR_+$ and the objective functional $\E_T(u_0) := \E(u(T))$,
find
\begin{equation*}
\tuuET  =  \mathop{\arg\max}_{u_0 \in \Xi_{\E_0}} \, \E_T(u_0), \qquad \text{where} \qquad
\Xi_{\E_0}  :=  \left\{u_0\in H^1(\Omega)\,\colon \; \int_0^1 u_0 \, dx = 0, \ \E(u_0) = \E_0 \right\}.
\end{equation*} 
\end{problem}
\noindent
The idea behind this problem is to find optimal initial data $\tuuET$
with prescribed enstrophy $\E_0$ that at the given time $T$ produces
the largest enstrophy $\E_T(\tuuET)$. We emphasize that in involving
the flow evolution on $[0,T]$, Problem \ref{pb:maxET1D} is
fundamentally different, and arguably harder to solve, than Problem
\ref{pb:maxdEdt1D} where the instantaneous only amplification of
enstrophy is considered. Problem \ref{pb:maxET1D} was solved in
\cite{ap11a} with $\nu = 10^{-3}$ for a broad range of values of
$\E_0$ and $T$ using the adjoint-based gradient-ascent method
described in Appendix \ref{sec:opt} and some results are summarized in
Figures \ref{fig:maxET1D}a and \ref{fig:maxET1D}b. As is evident
from Figure \ref{fig:maxET1D}a, the optimal initial data $\tuuET$
obtained for a fixed enstrophy $\E_0$ and a short time window $T$
features a steep front and hence resembles the instantaneous
maximizers $\tuuE$ found in \cite{ld08} by solving Problem
\ref{pb:maxdEdt1D}, however, as $T$ increases it gradually turns into
a rarefaction wave. By maximizing the results presented in Figure
\ref{fig:maxET1D}b with respect to $T$ at fixed values of $\E_0$, we
obtain the relation
\begin{equation}
\max_T \E_T(\tuuET)  \sim 11.488 \, \E_0^{1.531} \quad \text{as} \quad \E_0 \rightarrow \infty,
\label{eq:maxETap}
\end{equation}
where the exponent of $\E_0$ is lower, roughly by a factor of 2, than
the exponent 3 in the finite-time estimate \eqref{eq:Etbound1D}.  This
indicates that this estimate may not {be} sharp and could
possibly be improved by lowering the exponent of $\E_0$. We will
return to this question in Section \ref{sec:bound}.

Properties of extreme Burgers flows corresponding to {the
  initial data $\tuuET$ obtained as local maximizers of Problem
  \ref{pb:maxET1D} were analyzed} by Pelinovsky \cite{p12,p12b}. In
particular, subject to {the additional assumption that the
  initial condition be given in terms of an odd $C^3$ function, an
  $\O(\E_0^{3/2})$ estimate} was established in \cite{p12} on the
maximum growth of enstrophy $\max_t \E(t)$, cf.~\eqref{eq:maxETap}. It
was obtained applying Laplace's method to produce an asymptotic
representation for large $\E_0$ of the solution to \eqref{eq:Burgers}
given in terms of the Cole-Hopf formula.  {These results provide
  a rigorous and quantitative justification for the behavior of
  Burgers flows with initial data $\tuuET$ obtained as local
  maximizers of Problem \ref{pb:maxET1D}.}

Since Problem \ref{pb:maxET1D} is nonconvex and the numerical approach
employed to solve it relies on local optimality conditions,
cf.~Appendix \ref{sec:opt}, we of course cannot guarantee that the
solutions found for any $\E_0$ and $T$, cf.~Figures
\ref{fig:maxET1D}a,b, are global maximizers. However, the results
reported in \cite{ap11a} were obtained following a thorough search
involving the use of many different, mutually orthogonal (in the
function space $H^1(\Omega)$), and random initial guesses $u^{(0)}$.
The optimal initial conditions shown in Figure \ref{fig:maxET1D}a are
in fact nonunique maximizers as their rescaled copies $(1/m) \tuuET(m
x)$, $x\in [0,1]$, $m=2,3,\dots$, were also found to be local
maximizers, but characterized by smaller values of $\E_T$.
{Further support for the conjecture that at least up to a
  certain value of $\E_0$ the maximizers presented in Figure
  \ref{fig:maxET1D}a are in fact global maximizers was provided in
  \cite{FantuzziGoluskin2020} where upper bounds on $\max_{t \ge 0}
  \E(t)$ revealing behavior consistent with \eqref{eq:maxETap} were
  obtained based on a Galerkin truncation of the Burgers system
  \eqref{eq:Burgers}. We will discuss this important point in more
  detail in Section \ref{sec:bound}.}

An intriguing question, originally raised by Flandoli \cite{f15}, is
how extreme or singular behavior possible in solutions of hydrodynamic
models may be affected by stochastic forcing. More specifically, the
question is whether via some interaction with the nonlinearity {and
  dissipation present in the system} such stochastic forcing may
enhance or weaken the growth of certain solution norms as compared to
the deterministic case. In particular, in the case of systems
exhibiting finite-time blow-up in the deterministic setting it is
interesting to know whether noise may accelerate or delay the
formation of a singularity, or perhaps even prevent it entirely
\cite{f15}.  The question how colored additive noise in 1D Burgers
equation affects the dependence of the maximum attained enstrophy
$\max_{t \ge 0} \E(t)$ on $\E_0$ was investigated using stochastic
Monte-Carlo techniques in \cite{PocasProtas2018}. It was shown however
that the expected values of the enstrophy in stochastic Burgers flows
with the optimal initial conditions $\tuuET$ exhibit the same
power-law dependence on the initial enstrophy $\E_0$ as in the
deterministic case, cf.~\eqref{eq:maxETap}.

An interesting generalization of system \eqref{eq:Burgers} is the
fractional Burgers system
\begin{subequations}
\label{FBE}
\begin{align}
\partial_{t}u+\frac{1}{2}\,\partial_{x}u^{2}+\nu\,(-\Delta)^{\alpha} u=0 \quad & \mbox{in}\ (0,T]\times \Omega,\label{FBEa}\\
u(0) = u_0&,\label{FBEc}
\end{align}
\end{subequations}
where $(-\Delta)^{\alpha}$, $\alpha\in[0,1]$, is the fractional
Laplacian defined for sufficiently smooth functions $v \; : \; \Omega
\rightarrow \RR$ in terms of the relation $\left[
  \widehat{(-\Delta)^\alpha v} \right]_k := |k|^{2\alpha} [
\widehat{v} ]_k$, $k \in \ZZ$, in which $[\widehat{v} ]_k$ is the
Fourier coefficient of $v$ with wavenumber $k$. As shown in
\cite{kns08}, system \eqref{FBE} admits globally-defined smooth
classical solutions in the subcritical ($\alpha \in (1/2,1]$) and in
the critical ($\alpha = 1/2$) regime. On the other hand, finite-time
blow-up occurs in the supercritical regime ($\alpha \in [0,1/2)$). The
fractional Burgers system is thus a useful simple model to study
singular behavior, especially given the fact that the 3D Navier-Stokes
system is also known to be globally well posed in the classical sense
in the presence of fractional dissipation with exponents $\alpha \ge
5/4$ \cite{kp12}. The fractional Burgers system \eqref{FBE} has also
been studied in connection with turbulence \cite{Boritchev2018}.

Generalizations of the instantaneous estimate \eqref{eq:dEdtbound1D}
for the case of the fractional Burgers system \eqref{FBE} have been
obtained in \cite{Yun2018}. It was shown that the dependence of the
bounds on the enstrophy rate of growth $d\E/dt$ on $\E_0$ has the same
global form $\sigma \E_0^\gamma$ in the subcritical, critical and
parts of the supercritical regime with the exponent $\gamma$
increasing without bound as the fractional dissipation exponent
$\alpha$ is reduced from 1 (where $\gamma = 5/3$,
cf.~\eqref{eq:dEdtbound1D}) to 1/4.  Moreover, by solving numerically
a variant of Problem \ref{pb:maxdEdt1D}, these new estimates were
shown to be sharp (up to numerical prefactors). Finally, singularity
formation in the supercritical regime and transient behavior in the
subcritical case were studied numerically using Monte-Carlo methods in
fractional Burgers flows subject to additive colored noise in
\cite{RamirezProtas2021}.  The main finding was that there was no
evidence for the noise to regularize the evolution by suppressing
blow-up in the supercritical regime, or for the noise to trigger
blow-up in the subcritical regime.  However, as the noise amplitude
becomes large, the blow-up times in the supercritical regime
(understood as a random variable) were shown to exhibit an
increasingly non-Gaussian behavior.

\section{Estimates for the 2D Navier-Stokes System}
\label{sec:2D}

Denoting $\omega(t,\x) := \bomega(t,\x)\cdot \e_3$ the vorticity
component perpendicular to the plane of motion, where $\e_3$ is the
corresponding unit vector of the Cartesian coordinate system, the 2D
Navier-Stokes system can be written as
\begin{subequations}
\label{eq:NS2D}
\begin{alignat}{2}
 \Dpartial{\omega}{t} + J(\omega,\psi) & = \nu \Delta \omega &  \qquad &
\textrm{in} \ (0,T] \times \Omega, \label{eq:NS2Da} \\
 - \Delta \psi & = \omega & \quad &
\textrm{in} \ (0,T] \times \Omega, \label{eq:NS2Db} \\
 \omega(0) & = \omega_0 & \quad & \label{eq:NS2Dc}
\end{alignat}
\end{subequations}
where $\Omega = \TT^2$ is a doubly-periodic domain, $\psi$ the
streamfunction, whereas $J(f,g) := \partial_x f \, \partial_y g -
\partial_y f \,
\partial_x g$ defined for some functions $f,g \; : \; \Omega
\rightarrow \RR$ is the Jacobian determinant. As is well known
\cite{kl04}, system \eqref{eq:NS2D} is globally well posed in the
classical sense.

In the absence of vortex stretching in \eqref{eq:NS2Da}, the cubic term
responsible for the production of enstrophy in \eqref{eq:dEdt}
vanishes identically, such that for 2D flows on domains without solid
boundaries we have $d\E(\psi(t))/dt \le 0$, $t \ge 0$ (for
convenience, here we assume the streamfunction $\psi$ to be the main
state variable). Thus, in such cases the enstrophy is a nonincreasing
function of time and hence is rather uninteresting.

On the other hand, by computing the gradient of equation
\eqref{eq:NS2Da} we  obtain the equation describing the evolution of the
vorticity gradient $\bnabla\omega$
\begin{equation}
  \Dpartial{\bnabla\omega}{t} + (\u\cdot\bnabla)\bnabla\omega
  = \nu \Delta \bnabla\omega - \left[\bnabla\u\right]^T \cdot \bnabla\omega,
\label{eq:gradw}
\end{equation}
where the velocity field is given by $\u = \bnabla^{\perp} \psi$ with
$\bnabla^{\perp} := \left[ \partial / \partial_y, - \partial /
  \partial_x \right]$ and the palinstrophy
\begin{equation}
\P(\psi(t)) := \frac{1}{2}\int_{\Omega} |\bnabla \Delta \psi(t,\x)|^2 \, d\Omega
\label{eq:P}
\end{equation}
plays the role of ``energy''. Since equation \eqref{eq:gradw} features
a quadratic stretching term $\left[\bnabla\u\right]^T \cdot
\bnabla\omega$, palinstrophy may exhibit nontrivial growth in 2D
Navier-Stokes flows, as opposed to energy and enstrophy.  Hence, it
serves as a key measure of extreme behavior possible in such flows and
its rate of growth describing the build-up of vorticity gradients can
be obtained from \eqref{eq:gradw} as
\begin{equation}
\frac{d\P(\psi(t))}{dt} = \int_{\Omega} J(\Delta\psi,\psi) \Delta^2 \psi\, d\Omega
- \nu \, \int_{\Omega} (\Delta^2 \psi)^2 \, d\Omega =: \R_{\P}(\psi).
\label{eq:RP}
\end{equation}

In analogy with the results discussed in Section \ref{sec:1D}, the
goal of this study was to characterize the largest growth of
palinstrophy possible instantaneously and in finite time.  As a first
step, the following estimate on the rate of growth of palinstrophy was
obtained in \cite{ap13a}
\begin{equation}
\frac{d\P}{dt}  \le \frac{C}{\nu} \,\K^{\frac{1}{2}}\, \P^{\frac{3}{2}}.
\label{eq:dPdt}
\end{equation}
We note that, in contrast to the estimates on the rate of growth of
enstrophy in 1D and in 3D, cf.~\eqref{eq:dEdt_estimate_E} and
\eqref{eq:dEdtbound1D}, the upper bound in \eqref{eq:dPdt} is a
function of two quantities, i.e., the energy $\K$ and palinstrophy
$\P$. The former quantity could be eliminated from \eqref{eq:dPdt} in
favor of $\P$ using nested Poincar\'e's inequalities $\K \le (2
\pi)^{-4} \P$ which would give $d\P / dt \le (C /\nu) \P^2$, however,
sharpness would be lost in this process. Estimate \eqref{eq:dPdt} was
refined in \cite{ayala_doering_simon_2018} where a sharper form of the
prefactor dependent on $\K$ was obtained
\begin{equation}
\frac{d\P}{dt}  \le \left(a + b\sqrt{\ln Re +c} \right) \, \P^{\frac{3}{2}} \quad \text{with}
\quad  a = 0, \ b = \sqrt{2\pi}, \ c = - \ln\left(\frac{2}{\sqrt{\pi}}\right)
\label{eq:dPdt2}
\end{equation}
and with the Reynolds number defined as $Re := \K^{1/2} / \nu$.  We
add that, as was shown in \cite{td06} (see also \cite{ap13a}), some
other estimates on $d\P / dt$ can be obtained, but they involve bounds
on quantities such as $\|\Delta\omega\|_{L^2(\Omega)}$ and
$\|\omega\|_{L^\infty(\Omega)}$ which are hard to control. Estimates
on the rate of growth of palinstrophy in the presence of external body
forces were obtained in \cite{dfj10}. 

\begin{figure}[t]
\begin{center}
\mbox{
\subfigure[]{\includegraphics[width=0.5\textwidth]{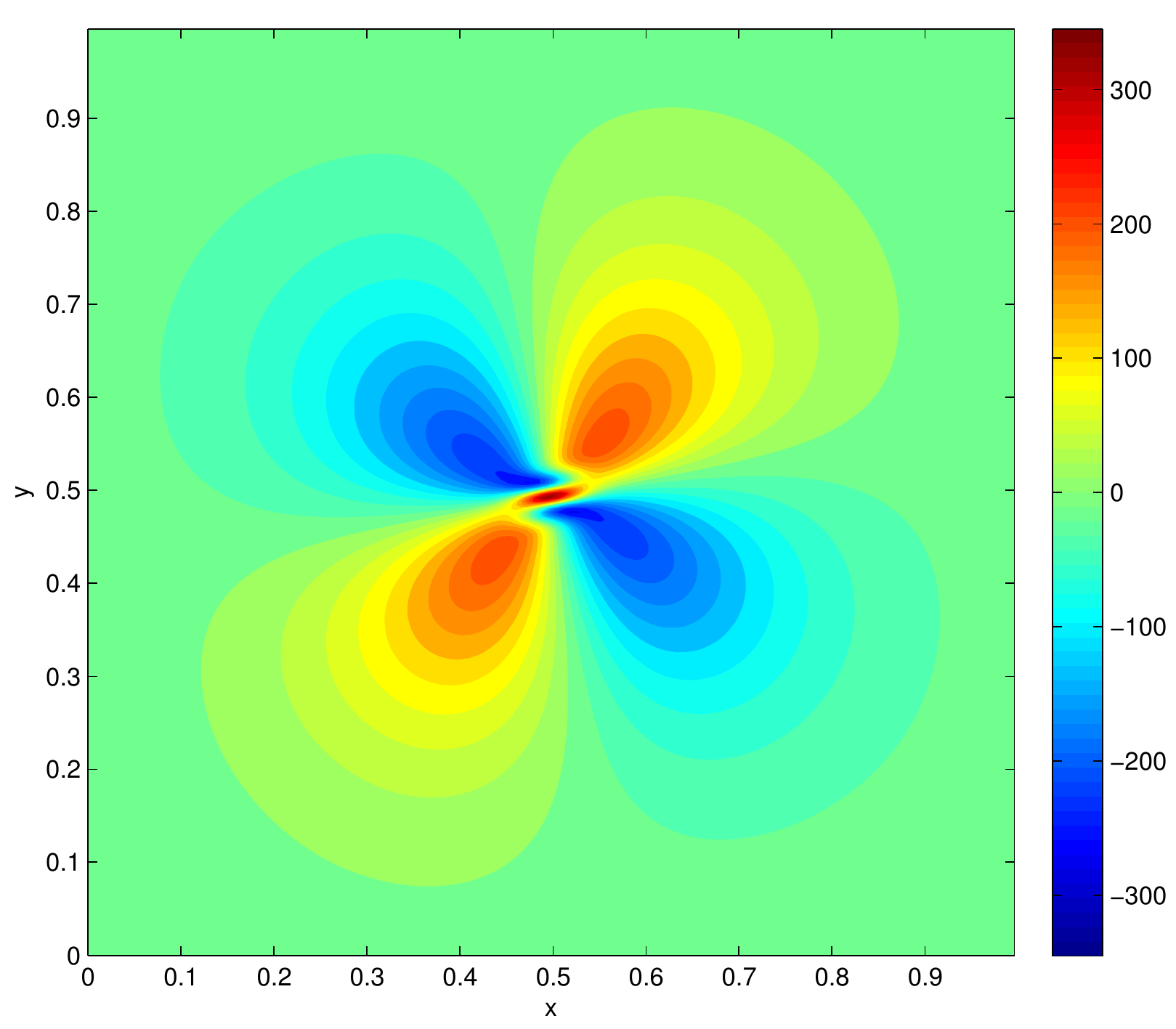}}\quad
\subfigure[]{\includegraphics[trim=0 1cm 0 -1cm, width=0.475\textwidth]{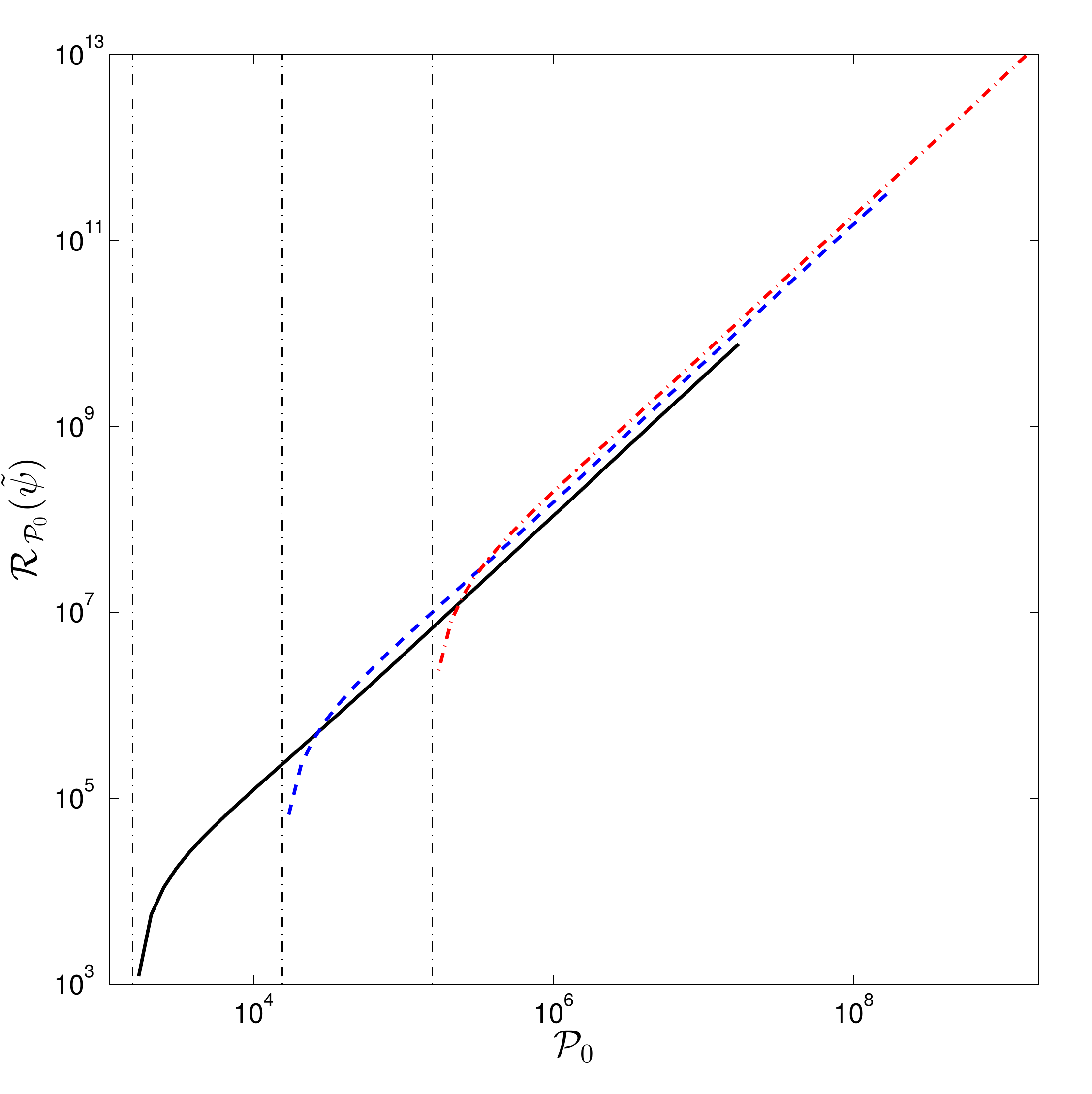}}}
\caption{(a) Vorticity field $-\Delta\tpsiKP$ solving Problem
  \ref{pb:maxdPdt} with $\K_0 = 10$ and $\P_0 =1.5585\cdot 10^6$. (b)
  Dependence of the maximum palinstrophy rate of growth
  $\R_{\P_0}(\tpsiKP)$ on $\P_0$ for $\K_0 = 10^0, 10^1$ and $10^2$
  with the vertical lines representing the corresponding Poincar\'e
  limits $(2 \pi)^{4} \K_0$. }
\label{fig:maxdPdt}
\vspace*{-0.8cm}
\end{center}
\end{figure}

By integrating the instantaneous estimate \eqref{eq:dPdt2} with
respect to time, the following finite-time bound was obtained in
\cite{ayala_doering_simon_2018}
\begin{equation}
\max_{t \ge 0} \P(\psi(t)) \le \Phi(Re_0) \P_0 \qquad \text{with}
\qquad   \Phi(Re_0) :=  \left(1 + \frac{a + b\sqrt{\ln Re_0 +c}}{4} Re_0 \right)^2,
\label{eq:Pt2}
\end{equation}
where $Re_0 := \K_0^{1/2} / \nu$ and $\P_0 := \P(\psi(0))$. In order
to assess sharpness of instantaneous estimates
\eqref{eq:dPdt}--\eqref{eq:dPdt2} with respect to $\K$ and $\P$, the
following optimization problem was formulated in \cite{ap13a}
\begin{problem}
\label{pb:maxdPdt}
Given $\K_0,\P_0\in\RR_+$ and the objective functional
$\R_{\P_0}(\psi)$, cf.~\eqref{eq:RP}, find
\begin{align*}
\tpsiKP & = \mathop{\arg\max}_{\psi\in\W_{\K_0,\P_0}} \,\R_{\P_0}(\psi), \qquad \text{where}  \\ 
\W_{\K_0,\P_0} & = \left\{\psi \in H^4(\Omega) : \frac{1}{2}\int_\Omega|\bnabla\psi|^2\,d\Omega = \K_0, \ \frac{1}{2}\int_\Omega|\bnabla\Delta\psi|^2\,d\Omega = \P_0 \right\}.
\end{align*}
\end{problem}
\noindent 
We emphasize that in contrast to Problem \ref{pb:maxdEdt1D}, Problem
\ref{pb:maxdPdt} involves two constraints which {is motivated by
  the structure of the upper bounds in estimates
  \eqref{eq:dPdt}--\eqref{eq:dPdt2} and} makes it harder to solve
numerically, cf.~comments at the end of Appendix \ref{sec:optP4}.
{Local maximizers of Problem \ref{pb:maxdPdt} with $\nu =
  10^{-3}$ were found in \cite{ap13a}} for a broad range of values of
$\K_0$ and $\P_0$, where we focused on the dependence of
$\R_{\P_0}(\tpsiKP)$ on $\P_0$ with the kinetic energy $\K_0$ held
fixed. A representative maximizer $\tpsiKP$ is shown in terms of the
corresponding vorticity field $-\Delta \tpsiKP$ in Figure
\ref{fig:maxdPdt}a. As is evident from this figure, the optimal state
involves a quadrupole vortex generating a straining field that
stretches a vortex filament located at the center. As demonstrated in
\cite{ayala_doering_simon_2018}, for fixed $\K_0$ the optimizers
$\tpsiKP$ are self-similar with respect to $\P_0$, i.e., they admit
the representation $\tpsiKP = \P_0^{\beta} \Psi(\P_0^q \x)$, where the
rational exponents $\beta$ and $q$ are determined from the constraints
and the optimality conditions in Problem \ref{pb:maxdPdt}, whereas
$\Psi$ is a function independent of $\P_0$, but depending on $\K_0$.
We also note that in the small-palinstrophy limit defined by
Poincar\'e's inequality $\P_0 \rightarrow (2 \pi)^{4} \K_0$, the cubic
term in \eqref{eq:RP} vanishes which simplifies Problem
\ref{pb:maxdPdt} since the objective function becomes quadratic. This
limiting problem can be solved in closed form using the method of
Lagrange multipliers with the maximizers $\tpsiKP$ having the form of
eigenfunctions of the Laplacian.

In \cite{ap13a} we were interested in sharpness of the instantaneous
estimate \eqref{eq:dPdt} with respect to $\P_0$ and this is verified
in Figure \ref{fig:maxdPdt}b where we plot $\R_{\P_0}(\tpsiKP)$ as
function of $\P_0$ for different values of $\K_0$, revealing a
power-law dependence of the form $\R_{\P_0}(\tpsiKP) \sim
\O(\P_0^{3/2})$ for large $\P_0$. Sharpness of the refined estimate
\eqref{eq:dPdt2} with respect to the second parameter $\K_0$ was then
established in \cite{ayala_doering_simon_2018} by analyzing solutions
of Problem \ref{pb:maxdPdt}, although the parameters $a$, $b$ and $c$
were found to have values different from the values given in
\eqref{eq:dPdt2}.

As regards the companion question about sharpness of the finite-time
estimate \eqref{eq:Pt2}, in \cite{ayala_doering_simon_2018} it was
shown that flow evolutions corresponding to the optimal initial data
$\tpsiKP$ obtained by solving Problem \ref{pb:maxdPdt} for fixed
$\K_0$ and different $\P_0$ saturate the bound in estimate
\eqref{eq:Pt2} with respect to $\P_0$, in the sense that the maximum
attained values of the palinstrophy $\max_{t \ge 0} \P(\psi(t))$ grow
in proportion to $\P_0$ when $\K_0$ remains fixed. On the other hand,
dependence of the prefactor $\Phi(Re_0)$ on $Re_0$, or on $\K_0$, was
found to be more nuanced, which can be attributed to the fact that the
flow trajectories considered correspond to initial data which is
optimal in the instantaneous sense only. The vortex-dynamics
mechanisms realizing the extreme flow behavior discussed above were
analyzed in \cite{ap13b}. As is evident from Movie 1, the stretching
of three thin parallel vortex filaments is the key effect responsible
for the build-up of the palinstrophy. We add that since $d\E / dt =
- \nu \P$ the question about the maximum growth of palinstrophy is
related to the problem of the enstrophy dissipation vanishing in the
limit $\nu \rightarrow 0$ in 2D turbulence
\cite{td06,JeongYoneda2021}.

Finally, we add that on bounded domains there are additional terms in
expression \eqref{eq:dEdt} for the rate of growth of enstrophy in the
form of integrals over the domain boundary $\partial\Omega$, such that
in 2D Navier-Stokes flows on such domains the enstrophy can grow.
While we are unaware of any a priori estimates on $d\E/dt$ on bounded
domains in 2D, the extreme behavior of this quantity can be studied by
solving suitable optimization problems and some preliminary results in
this direction were reported in \cite{Sliwiak2017}.

\section{Estimates for the 3D Navier-Stokes System}
\label{sec:3D}

The question about sharpness of the instantaneous estimate
\eqref{eq:dEdt_estimate_E} was considered by Lu \& Doering in
\cite{ld08} who formulated and {studied} the following
optimization problem
\begin{problem}
\label{pb:maxdEdt3D}
  Given $\E_0\in\RR_+$ and the objective functional $\R_\E(\u)$, cf.~\eqref{eq:dEdt},
find
\begin{equation*}
\tuE  =  \mathop{\arg\max}_{\u \in \mathcal{S}_{\E_0}} \, \R_\E(\u), \qquad \text{where} \qquad
\mathcal{S}_{\E_0}  :=  \left\{\u\in H^2(\Omega)\,\colon \; \bnabla\cdot\u =0, \ \int_{\Omega} \u \, d\x = \0, \ \E(\u) = \E_0 \right\}.
\end{equation*} 
\end{problem}
\noindent
We remark that the numerical approach adopted in \cite{ld08} was
somewhat different from the methodology described in Appendix
\ref{sec:opt} in that it relied on a ``discretize-then-optimize''
formulation wherein Problem \ref{pb:maxdEdt3D} was first discretized
with a Fourier-Galerkin method which then lead to an optimization
problem in a finite dimension. {Using this approach over a range
  of values of $\E_0$ and with $\nu = 10^{-2}$, Lu \& Doering found
  two branches of locally maximizing solutions $\tuE$ of Problem
  \ref{pb:maxdEdt3D} , with one branch} characterized by the
relation\footnote{{Due to the presence of the factor 1/2 in
    \eqref{eq:E}, the coefficient in \eqref{eq:maxR3D} differs from
    that given in \cite{ld08}.}}
\begin{equation}
\R(\tuE) \sim {3.59 \cdot 10^{-3}} \E_0^{2.997} \quad \text{as} \quad \E_0 \rightarrow \infty.
\label{eq:maxR3D}
\end{equation}
Thus, the maximizers on this branch, which interestingly have the form
of two colliding nearly axisymmetric vortex rings, saturate estimate
\eqref{eq:dEdt_estimate_E} in the sense that the rate at which these
maximizers produce enstrophy increases in proportion to $\E_0^3$
(although the numerical prefactor in \eqref{eq:maxR3D} is smaller than
the one in estimate \eqref{eq:dEdt_estimate_E} by about 9 orders of
magnitude).  These maximizers are strongly localized such that as the
enstrophy increases the characteristic radius of the vortex rings
vanishes as $\O(\E_0^{-1})$ \cite{ap16}. The maximizers associated
with the second branch were characterized by the asymptotic relation
$\R(\tuE) \sim {0.299} \, \E_0^{1.78}$ and involved vorticity
concentrated in four rod-like regions.

Problem \ref{pb:maxdEdt3D} was revisited in \cite{ap16} where we
{recomputed the asymptotically dominating branch with more
  accuracy which allowed us to slightly} improve the prefactor in
\eqref{eq:maxR3D} to $3.72\cdot 10^{-3}$. {We also considered
  the problem in the limit $\E_0 \rightarrow 0$ in which it} was shown
to admit closed-form solutions $\tuE$ in the form of divergence-free
eigenfunctions of the vector Laplacian.  One of these limiting
maximizers is the Taylor-Green vortex, which has been employed as the
initial data in a number of studies aimed at triggering singular
behaviour in both the Euler and Navier-Stokes systems
\cite{tg37,bmonmu83,b91,bb12}. It is interesting that the Taylor-Green
vortex arises as a solution of Problem \ref{pb:maxdEdt3D} in the limit
$\E_0 \rightarrow 0$.

\begin{figure}
\begin{center}
\mbox{\hspace*{-0.35cm}
\Bmp{0.5\textwidth}
\subfigure[]{\includegraphics[width=1.1\textwidth]{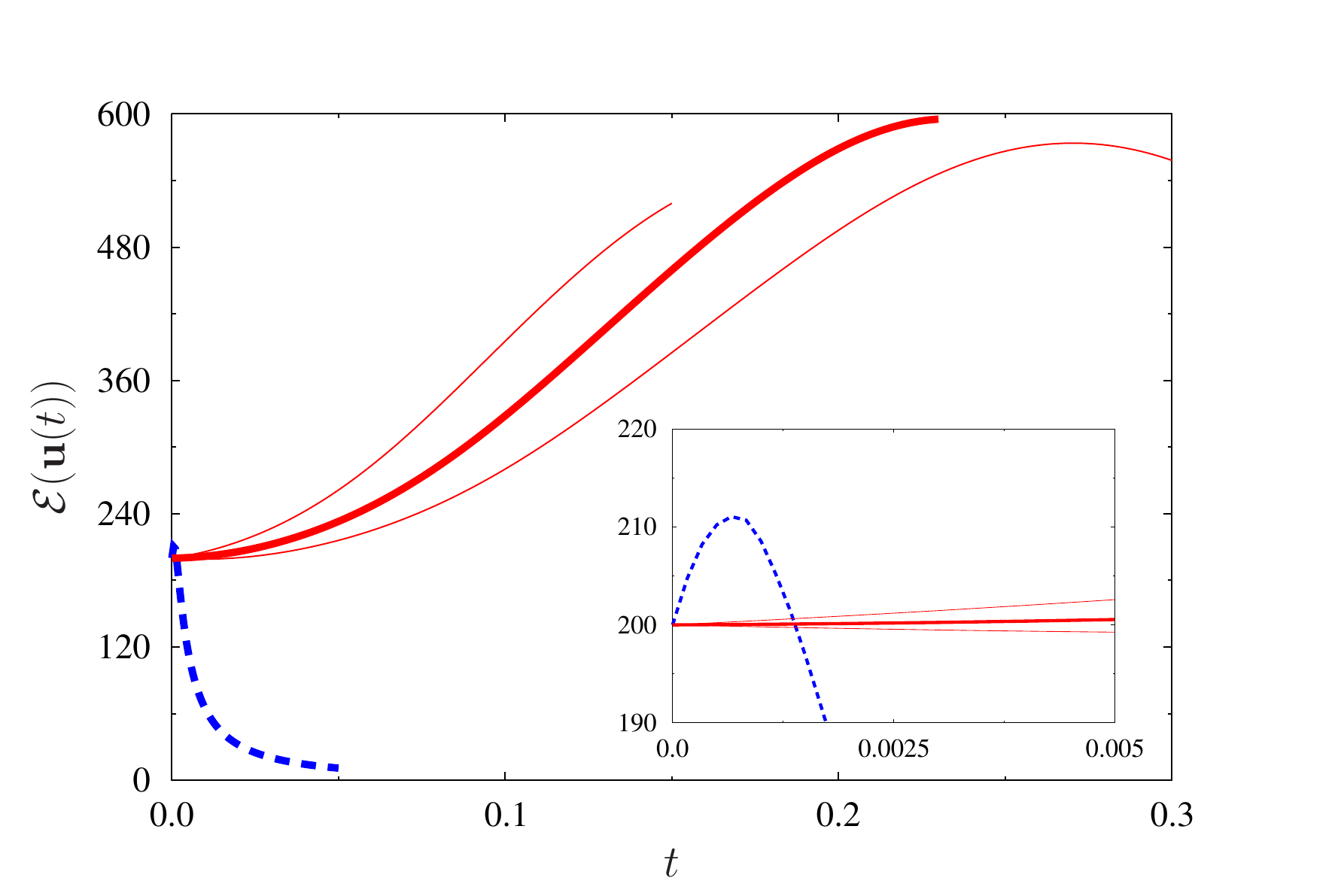}}
\Emp
\Bmp{0.5\textwidth}
\subfigure[]{\includegraphics[width=1.1\textwidth]{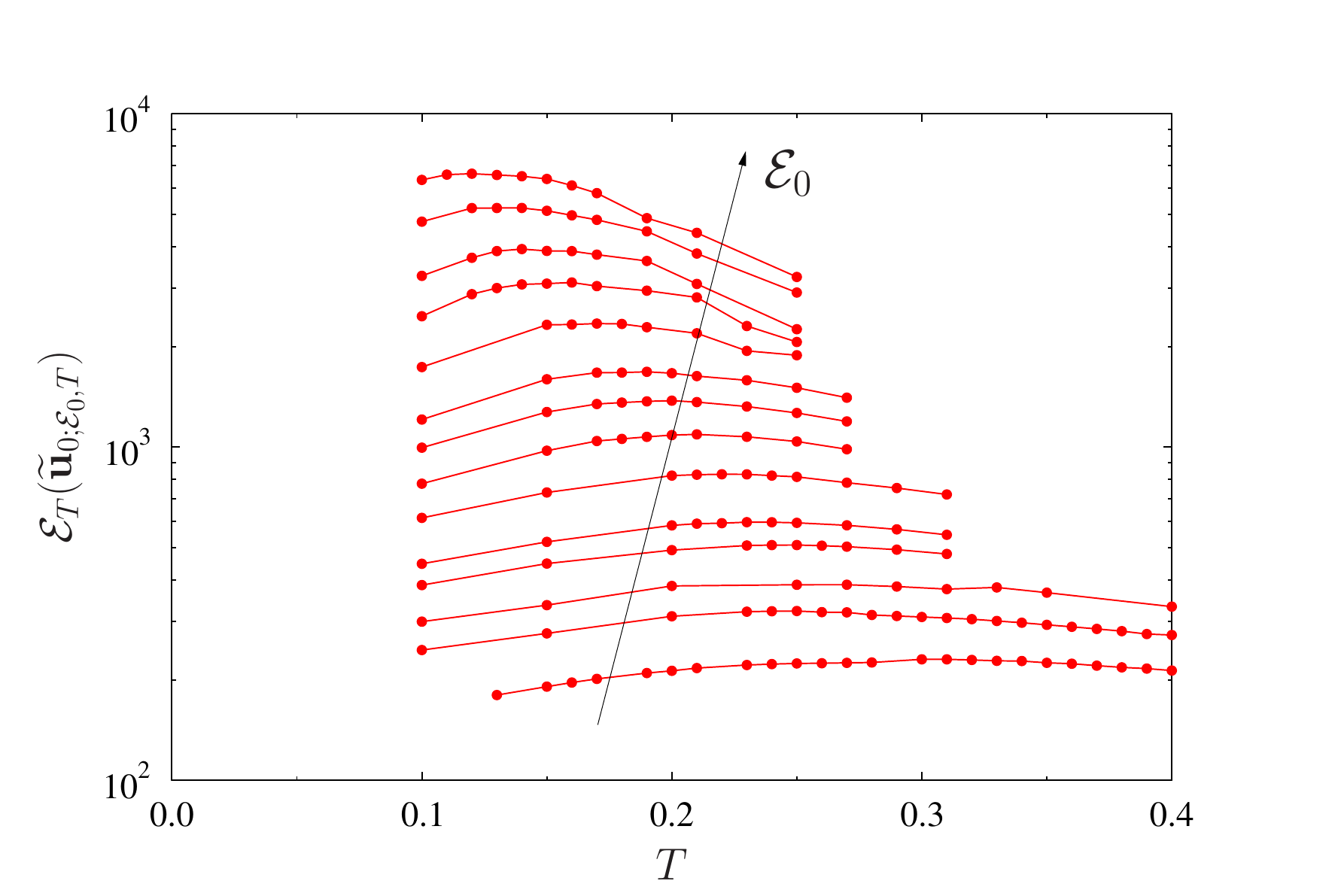}}
\Emp
}
\caption{(a) Enstrophy $\E(\u(t))$ as a function of time $t$ obtained
  from the solution of the Navier-Stokes system \eqref{eq:NS} with the
  initial condition $\u_0$ given by (blue dashed line) the maximizer
  $\tuE$ of Problem \ref{pb:maxdEdt3D} and (red solid lines) the
  asymmetric maximizers $\tuET$ of Problem \ref{pb:maxET3D} for $\E_0
  = 200$ and $T = 0.15,0.23,0.3$ (the curve corresponding to the
  optimal length of the time window $\tTE = 0.23$ is marked with a
  thick line whereas the inset represents magnification of the initial
  stages of evolution). (b) Maximum attained enstrophy $\E_T(\tuET)$
  as a function of the length $T$ of the window over which
  maximization is performed in Problem \ref{pb:maxET3D} for initial
  enstrophies $100 \le \E_0 \le 1000$ (the trend with the increase of
  $\E_0$ is indicated with an arrow).  Each curve corresponds to flow
  evolutions starting from asymmetric optimal initial conditions
  $\tuET$ with the same value of the initial enstrophy and different
  $T$ (solid symbols represent the values of $\E_0$ and $T$ for which
  {local maximizers of Problem \ref{pb:maxET3D} were found}).}
\label{fig:maxE3D}
\vspace*{-0.8cm}
\end{center}
\end{figure}

The time evolution of solutions of the Navier-Stokes system
\eqref{eq:NS} with the maximizers $\tuE$ of Problem \ref{pb:maxdEdt3D}
used as initial data was considered in \cite{ap16}. It was shown that
while at $t = 0$ the enstrophy in these flows is generated at the
maximum rate given in \eqref{eq:maxR3D}, this rate is very quickly
depleted such that in finite time only little enstrophy is produced,
cf~Figure \ref{fig:maxE3D}a. The flow evolution remains essentially
axisymmetric with the vortex rings approaching each other before
starting to diffuse. The key conclusion from these results is that if
a significant, let alone unbounded, growth of enstrophy is to be
achieved in finite time, it must be associated with initial data
$\u_0$ other than the extreme vortex states $\tuE$ saturating the
upper bound in estimate \eqref{eq:dEdt_estimate_E} on the
instantaneous rate of growth of enstrophy, cf.~\eqref{eq:maxR3D}.

More {specifically}, assuming the instantaneous rate of growth of
enstrophy in the form $d\E / dt = C \, \E^{\alpha}$ {with some
  prefactor} $C>0$, any exponent $\alpha > 2$ will cause $\E(\u(t))$
to become unbounded {at some} finite time {$t_0 = t_0(\alpha)$} if
this rate of growth is sustained over the interval $[0,t_0)$.  The
fact that there is no blow-up {when} $1 < \alpha \le 2$ follows from
the {observation} that one factor of $\E$ in
\eqref{eq:dEdt_estimate_E} can be bounded in terms of the initial
energy $\K_0$ using \eqref{eq:dKdt} as follows
\begin{equation}
\int_0^t \E(\u(s))\, ds = \frac{1}{2\nu} \left[ \K_0 - \K(\u(t))\right] \leq \frac{1}{2\nu} \K_0,
\label{eq:Kt}
\end{equation}
which upon applying Gr\"onwall's lemma to $d\E / dt = C \,
\E^{\alpha}$ with $\alpha = 2$ yields the bound
\begin{equation}
\max_{0 \le t \le T} \E(\u(t))  \le  \E_0\, \exp\left[C\int_0^T  \E(\u(s))\, ds\right] \le \E_0 \, \exp\left[\frac{C}{2\nu} \K_0\right].
\label{eq:Gronwall}
\end{equation}
Evidently, as the rate of growth of enstrophy slows down when $\alpha
\rightarrow 2^+$, {for blow-up to occur this minimum growth rate}
must be sustained over windows of time with increasing length, i.e.,
$t_0 \rightarrow \infty$ as $\alpha \rightarrow 2^+$. To assess the
feasibility of such a scenario, the following optimization problem was
considered in \cite{KangYunProtas2020}
\begin{problem}
\label{pb:maxET3D}
  Given $\E_0, T\in\RR_+$ and the objective functional $\E_T(\u_0) := \E(\u(T))$,
find
\begin{equation*}
\tuET  =  \mathop{\arg\max}_{\u_0 \in \Q_{\E_0}} \, \E_T(\u_0), \quad \text{where} \quad
\Q_{\E_0}  :=  \left\{\u_0\in H^1(\Omega)\,\colon  \bnabla\cdot\u_0 = 0, \ \int_{\Omega} \u_0 \, d\x = \0, \ \E(\u_0) = \E_0 \right\}.
\end{equation*} 
\end{problem}

\begin{figure}
\begin{center}
\mbox{\subfigure[$\tomega_1$]{\includegraphics[width=0.3\textwidth]{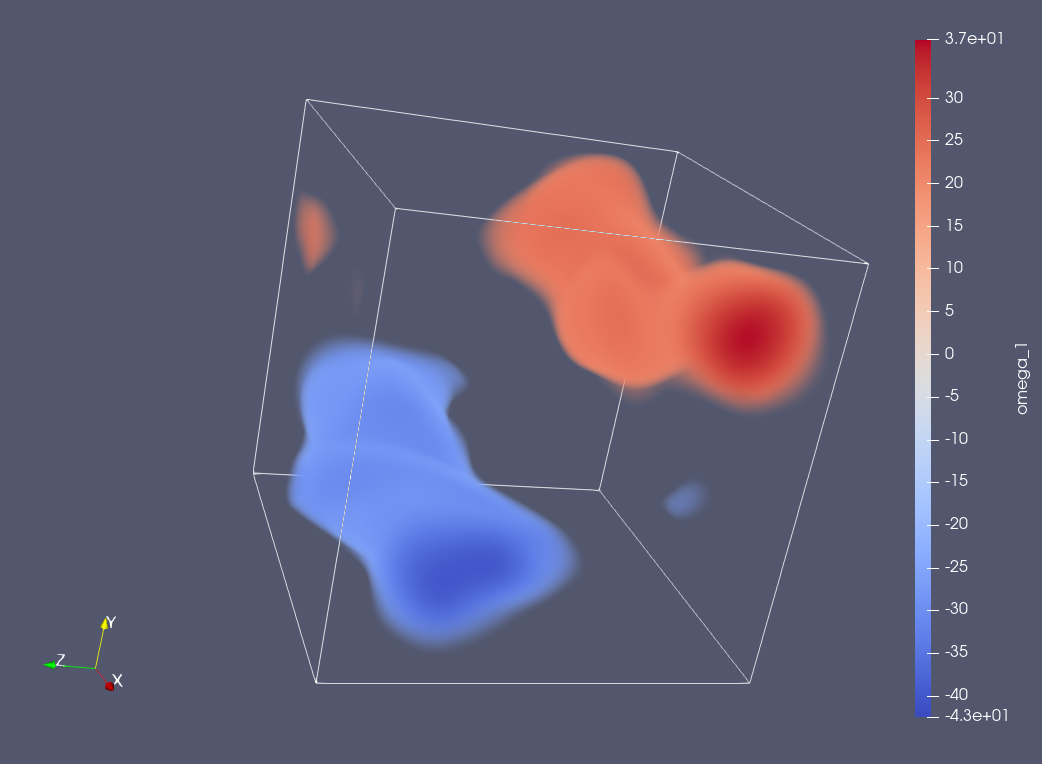}}\qquad
\subfigure[$\tomega_2$]{\includegraphics[width=0.3\textwidth]{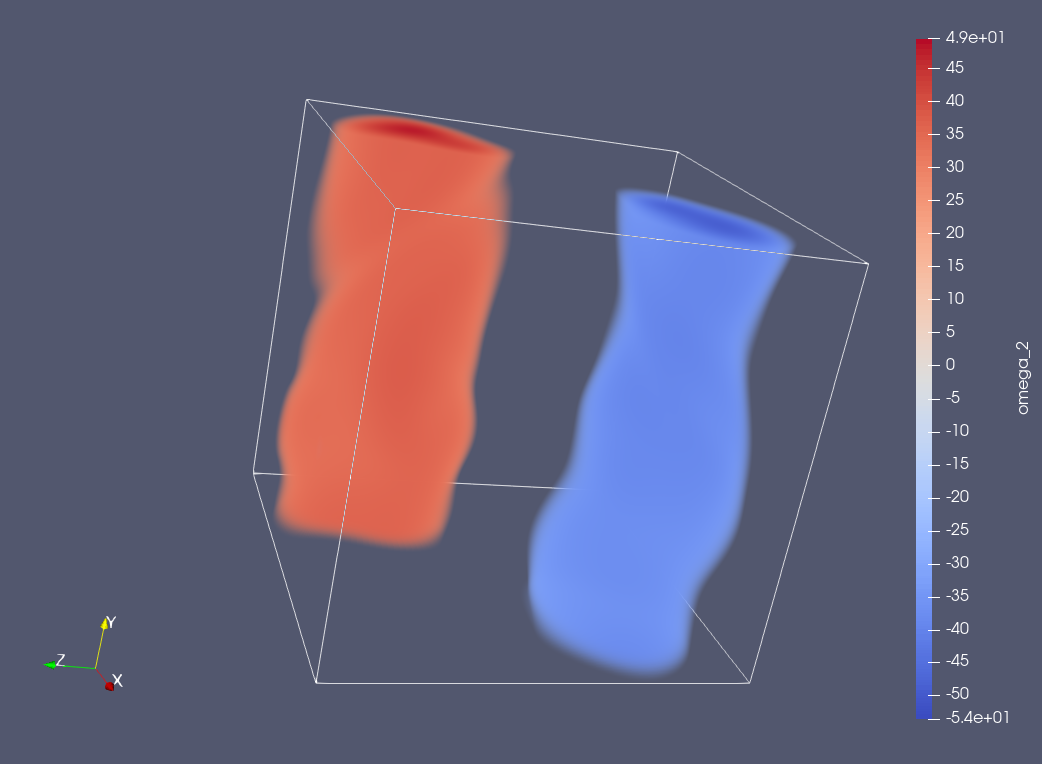}}\qquad
\subfigure[$\tomega_3$]{\includegraphics[width=0.3\textwidth]{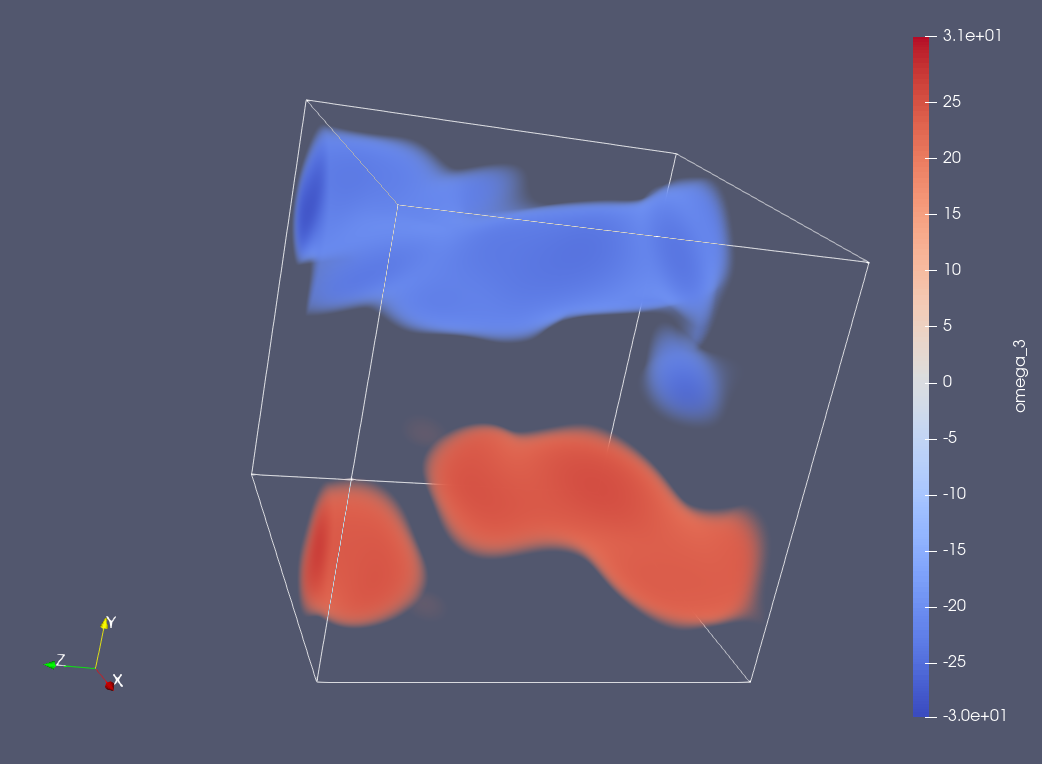}}}
\caption{Vorticity components of the {asymmetric} optimal initial
  condition $\tuEtT$ obtained by solving Problem \ref{pb:maxET3D} for
  the initial enstrophy $\E_0 = 500$ and the corresponding optimal
  length $\tTE = 0.17$ of the time interval. The time evolution of
  the flow corresponding to this initial condition is visualized in
  Movie 2.}
\label{fig:tubes}
\vspace*{-0.8cm}
\end{center}
\end{figure}

While Problem \ref{pb:maxET3D} is quite challenging from the
computational point of view, optimal solutions $\tuET$ were found in
\cite{KangYunProtas2020} for a range of values of $\E_0$ and $T$ with
$\nu = 10^{-2}$. They belong to two distinct branches, referred to as
``symmetric'' and ``asymmetric'', with the extreme flows corresponding
to the initial data $\tuET$ on the symmetric branch exhibiting
equipartition of enstrophy among the three Cartesian coordinate
dimensions. For large values of the initial enstrophy $\E_0$ the
asymmetric branch dominates in the sense that the corresponding
Navier-Stokes flows achieve higher values of $\E_T(\tuET)$ than the
flows with initial data on the symmetric branch for the same values of
$\E_0$ and $T$. The time evolution of the enstrophy $\E(\u(t))$ in the
extreme flows with the asymmetric initial conditions $\tuET$ obtained
for a fixed $\E_0 = 200$ and different time windows $T$ is shown in
Figure \ref{fig:maxE3D}a, where we see that in these flows a much
larger growth of enstrophy is achieved than in the flow with the
instantaneously optimal initial condition $\tuE$ obtained by solving
Problem \ref{pb:maxdEdt3D} for the same value of $\E_0$.
Interestingly, we notice that for some values of $T$ the enstrophy
$\E(\u(t))$ is in fact decreasing at early times before it starts to
grow. A typical asymmetric initial condition $\tuET$ is shown in
Figure \ref{fig:tubes} where it is evident that it has the form of
three perpendicular pairs of antiparallel vortex tubes. The evolution
of the flow corresponding to this initial condition is visualized in
Movie 2.  Interestingly, all extreme flows with initial conditions
found by solving Problem \ref{pb:maxET3D} have zero helicity
$\H(\u(t)) := \int_{\Omega} \u(t) \cdot \left(\bnabla \times \u(t)
\right) \, d\x = 0$, $t \ge 0$.

\begin{figure}
\begin{center}
\mbox{
\subfigure[]{\includegraphics[width=0.5\textwidth]{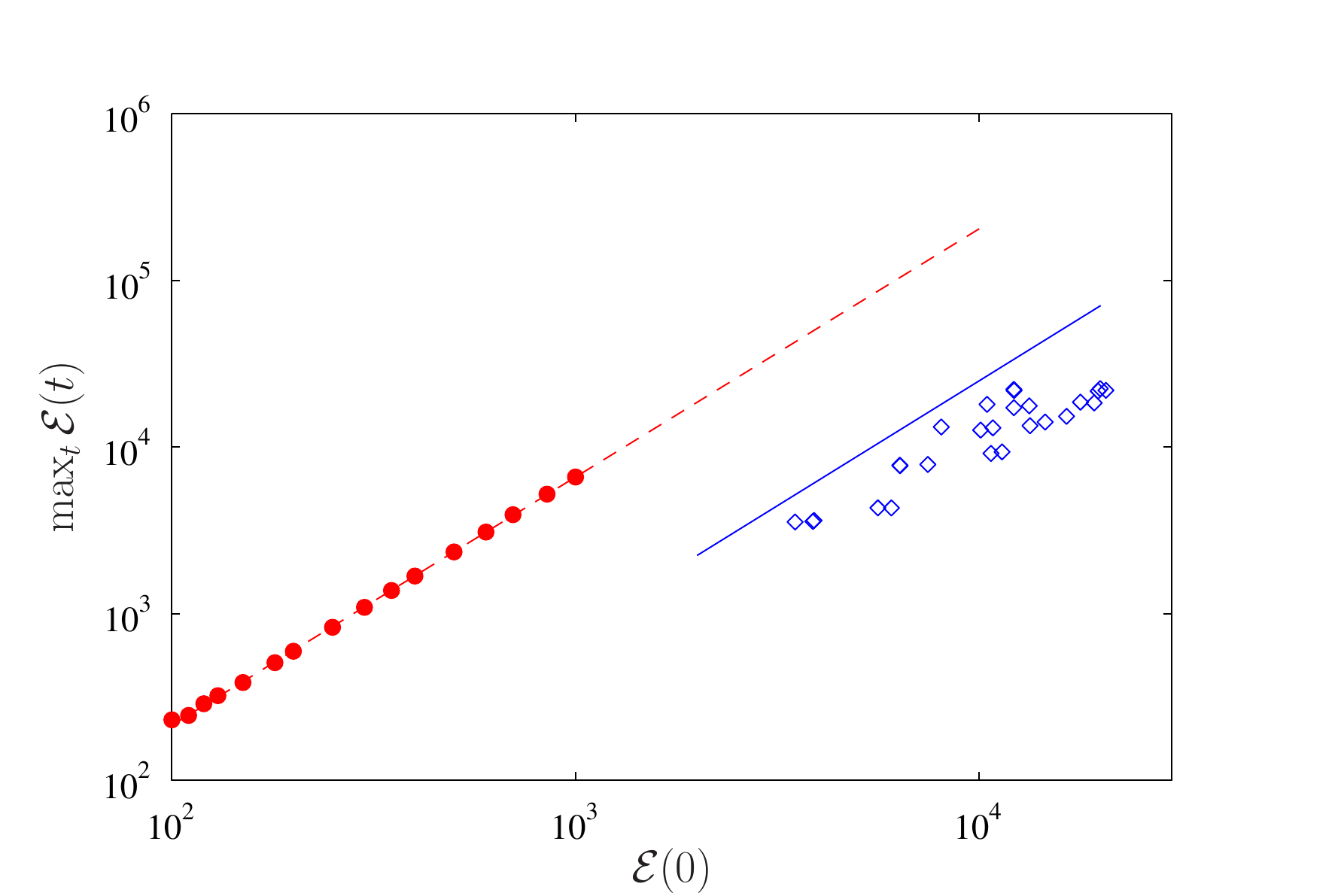}}\quad
\subfigure[]{\includegraphics[width=0.5\textwidth]{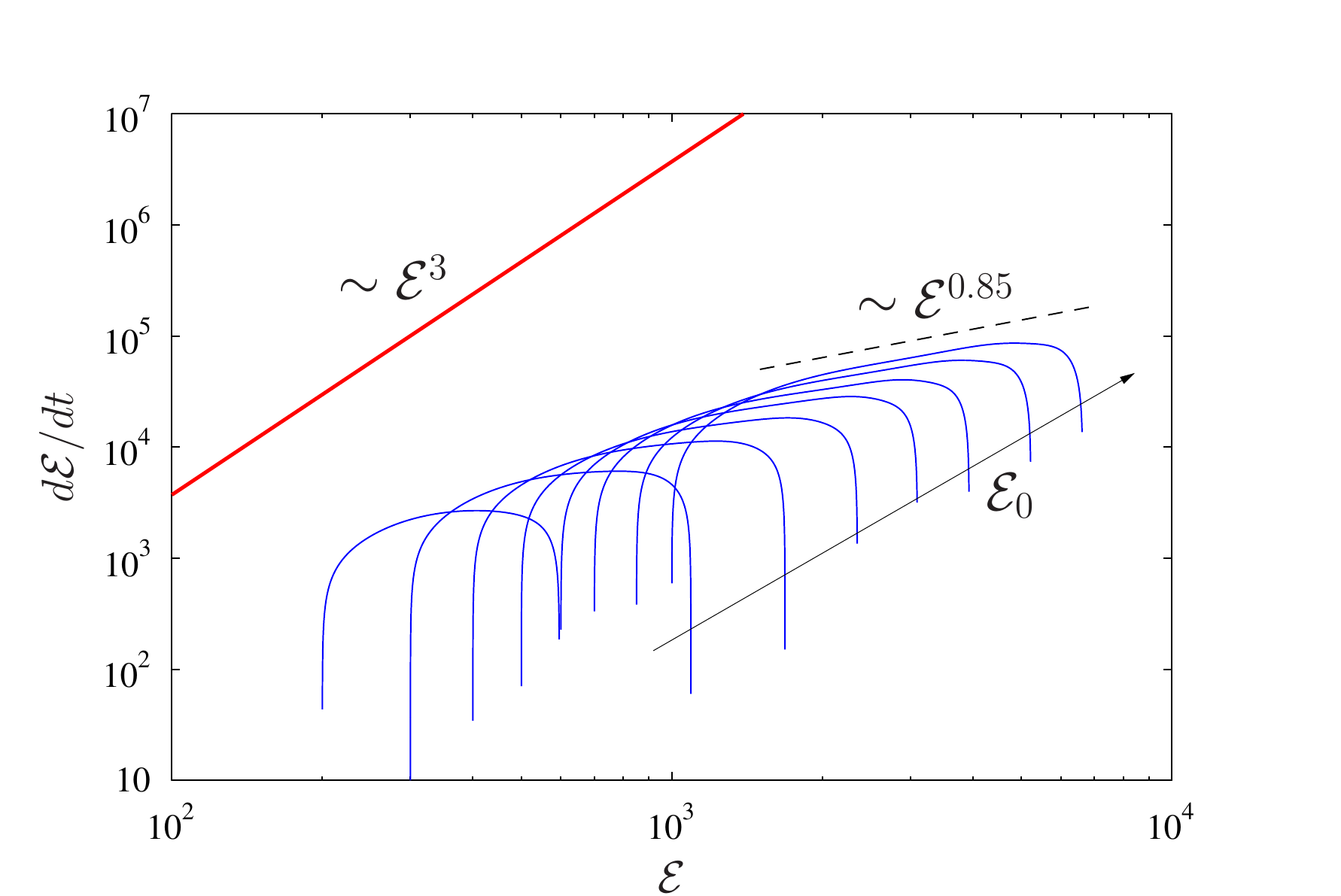}}
} 
\caption{(a) Dependence of the maximum attained enstrophy $\max_{t \ge
    0} \E(t)$ on the initial enstrophy $\E_0$ in Navier-Stokes flows
  with the optimal initial conditions (red solid circles) $\tuET$
  obtained by solving Problems \ref{pb:maxET3D} (the asymmetric
  branch) and (blue diamonds) $\tuBT$ obtained by solving Problem
  \ref{pb:PhiL4} (the partially-symmetric branch). {For the
    local maximizers of Problems \ref{pb:maxET3D}} each symbol
  corresponds to a different value of the constraint $\E_0$, and in
  all cases the results are presented for the optimization window with
  length $T$ producing the largest value of $\max_{t \ge 0} \E(t)$.
  {For Problem \ref{pb:PhiL4} the symbols correspond to local
    maximizers obtained with different values of $B$ and $T$, whereas
    the} straight lines represent the {relation} $\max_{t \ge
    0} \E(t) \approx C \E_0^{3/2}$ with different prefactors $C$. (b)
  Flow trajectories corresponding to the optimal initial data $\tuEtT$
  obtained by solving Problem \ref{pb:maxET3D} with different $\E_0
  \in [100,1000]$ shown using the coordinates $\{ \E, d\E/dt\}$ (blue
  solid lines with the arrow indicating the trend with the increase of
  $\E_0$). The thick red line represents the relation $d\E/dt =
  3.72\cdot 10^{-3} \, \E^3$ found in \cite{ap16}, whereas the dashed
  black line the relation $d\E/dt = 10^{2} \, \E^{0.85}$.}
\label{fig:maxETvsE0}
\vspace*{-0.8cm}
\end{center}
\end{figure}

No evidence has been found for unbounded growth of enstrophy that
would signal singularity formation, cf.~condition \eqref{eq:supE}, in
Navier-Stokes flows with initial data obtained as solutions of Problem
\ref{pb:maxET3D} for a broad range of values of $\E_0$ and $T$. The
maximum enstrophy values $\E_T(\tuET)$ attained with asymmetric
initial conditions $\tuET$ are shown as functions of the optimization
window $T$ for different values of $\E_0$ in Figure \ref{fig:maxE3D}b.
We see that the branches of maximizers corresponding to different
values of $\E_0$ all exhibit well-defined unique maxima attained at
times $\tTE = \argmax_{T>0} \E_T(\tuET)$, $\forall \E_0$, which
decrease as $\O(\E_0^{-1/2})$. These maximum enstrophy values
$\E_{\tTE}(\tuEtT)$ are plotted as a function of $\E_0$ in Figure
\ref{fig:maxETvsE0}a revealing a power-law relation
\begin{equation}
\max_{T > 0} \E_T(\tuET)  \ \sim \ \left( 0.224 \ \pm 0.006 \right) \, \E_0^{1.490 \, \pm 0.004}.
\label{eq:maxEt_vs_E0}
\end{equation}

In order to understand how close the flow evolutions corresponding to
the optimal initial data $\tuEtT$ come to saturating a priori bounds
on the rate of growth of enstrophy, cf.~\eqref{eq:dEdt_estimate_E}, in
Figure \ref{fig:maxETvsE0}b we plot the corresponding trajectories
using the coordinates $\{ \E, d\E/dt\}$, such that each trajectory is
parameterized by time $t$ (since the logarithmic scale is used,
initial parts of the trajectories when $d\E/dt \lessapprox 0$ are not
shown). The slope of the tangent to each of the curves thus represents
the exponent $\alpha$ characterizing the instantaneous rate of
enstrophy production $d\E/dt \sim \E^{\alpha}$. In Figure
\ref{fig:maxETvsE0}b we also indicate the relation $d\E/dt = 3.72\cdot
10^{-3} \, \E^3$ describing the maximum rate of enstrophy growth
realized by solutions of Problem \ref{pb:maxdEdt3D} \cite{ld08,ap16}.
We observe that the rate of growth of enstrophy achieved along the
trajectories corresponding to the optimal initial conditions $\tuEtT$
is at all times and for all values of $\E_0$ several orders of
magnitude smaller than the maximum rate of growth achieved by the
instantaneous maximizers $\tuE$. We also note that at the final stages
of the flow evolutions before the enstrophy maximum is reached at $t =
\tTE$ the enstrophy is amplified at an approximate rate $d\E/dt \sim
\E^{0.85}$, far below the minimum rate of growth $d\E/dt \sim
\E^{\alpha}$ with $2 < \alpha \le 3$ needed for enstrophy to become
unbounded in finite time.

In a recent study \cite{KangProtas2021} we have considered the
Ladyzhenskaya-Prodi-Serrin condition \eqref{eq:LPSblowup} focusing on
the case with $q = 4$ and $p = 8$, which is the pair of integer-valued
indices closest to the critical case with ${q} = 3$,
cf.~\eqref{eq:LPS3}. The goal was to search for potential finite-time
singularities in Navier-Stokes flows governed by \eqref{eq:NS} by
maximizing the quantity
\begin{equation}
\Phi_T(\u_0)  := \frac{1}{T} \int_0^T \| \u(\tau) \|_{L^4(\Omega)}^8 \, d\tau
\label{eq:Phi}
\end{equation}
with respect to the initial condition $\u_0$ where a natural function
space for the initial data would be $L^4(\Omega)$ and the constraint
would have the form $\| \u_0 \|_{L^4} = B$ for some $0 \le B \le
\infty$. However, from the computational point of view,
PDE-constrained optimization problems are formulated most conveniently
in a Hilbert space \cite{pbh04}. While there exist solution approaches
applicable in the more general setting of Banach spaces, e.g.,
\cite{protas2008}, they are significantly harder to use in practice.
Given the form of the constraint, we have chosen to formulate the
optimization problem in the ``largest'' Sobolev space with Hilbert
structure contained in $L^4(\Omega)$. From the Sobolev embedding
theorem in dimension 3 \cite{af05}, we deduce
\begin{equation}
H^s(\Omega) \hookrightarrow L^4(\Omega), \qquad s \ge \frac{3}{4},
\label{eq:HsL4}
\end{equation}
such that the largest Hilbert-Sobolev space embedded in $L^4(\Omega)$
is the space $H^{3/4}(\Omega)$. Thus, this leads to the following
optimization problem
\begin{problem}\label{pb:PhiL4}
  Given $B, T \in\mathbb{R}_+$ and the objective functional $\Phi_T(\u_0)$ from
  equation \eqref{eq:Phi}, find
\begin{align*}
\tuBT & =  \mathop{\arg\max}_{\u_0 \in \L_{B}} \, \Phi_T(\u_0), \quad \text{where} \\
\L_{B} & :=  \left\{\u_0\in H^{3/4}(\Omega)\,\colon\,\bnabla\cdot\u_0 = 0, \; \int_{\Omega} \u_0 \, d\x = \0,  \; \|\u_0\|_{L^4(\Omega)} = B \right\}.
\end{align*} 
\end{problem}

{Local maximizers of Problem \ref{pb:PhiL4} were found} in
\cite{KangProtas2021} for a range of values of $B$ and $T$ using a
numerical procedure which is an extension of the approach described in
Appendix \ref{sec:opt}, where the main modification concerned the
handling of the non-quadratic constraint in the definition of the
manifold $\L_B$. Two branches of maximizers were {discovered}
with partially symmetric and asymmetric optimal initial conditions
$\tuBT$. In this case as well no evidence was found for unbounded
growth of the quantity $\Phi_T(\tuBT)$ which would signal singularity
formation. The maximum enstrophy attained in the extreme flows with
the optimal initial conditions $\tuBT$ on the partially symmetric
branch obtained for different $B$ {and $T$} is plotted as
function of the initial enstrophy $\E_0$ in Figure
\ref{fig:maxETvsE0}a. It is intriguing to observe that {the
  envelope of these data points, obtained by maximizing the largest
  attained enstrophy over $B$ and $T$, is also described by the
  relation} $\max_{t>0} \E(t) \sim C \E_0^{3/2}$, i.e., the same as
found for flows corresponding to solutions of Problem
\ref{pb:maxET3D}, except that the prefactor $C$ is smaller than in
\eqref{eq:maxEt_vs_E0}.  Finally, by maximizing the quantity
$\Psi_T(\u_0) := \frac{1}{T} \int_0^T \| \u(\tau)
\|_{L^4(\Omega)}^{8/3} \, d\tau$ with respect to the initial data
$\u_0 \in H^{3/4}(\Omega)$ and subject to the constraint $(1/2)\| \u_0
\|_{L^2}^2 = \K_0$ for a range of $\K_0 > 0$ in an optimization
problem analogous to Problem \ref{pb:PhiL4} {we arrived at a
  conjecture that a priori estimate \eqref{eq:LPSbound} may} not sharp
and can possibly be improved by reducing the exponent of $\K_0$ in the
bound on the RHS. {The validity of this conjecture depends on
  whether the local maximizers of $\Psi_T(\u_0)$ we found are also
  global maximizers.}

\section{Relation to Bounding Approaches}
\label{sec:bound}

In this section we briefly discuss connections between the approaches
and results surveyed above and other techniques for quantifying the
extreme behavior possible in fluid flows. For brevity, we will assume
here the solution $\u(t) \in \X$, where $\X$ is a suitable
{Hilbert space of solutions} (finite or infinite dimensional),
satisfies the autonomous system $d{\u}(t)/dt = \f(\u(t))$ with some
$\f \; : \; \X \rightarrow \X$ and the initial condition $\u(0) = \u_0
\in X \subset \X$, where the set $X$ encodes the constraints imposed
on $\u_0$.  Denoting $\varphi \; : \; \X \rightarrow \RR$ the quantity
of interest in Problems \ref{pb:maxET1D}, \ref{pb:maxET3D} and
\ref{pb:PhiL4}, these problems can be expressed as
\begin{equation}
\overline{\varphi} := \sup_{\u_0 \in X} \varphi(\u(\cdot;\u_0)),
\label{eq:supphi}
\end{equation}
where $\u(t;\u_0)$ is the solution of the governing system at time $t$
corresponding to the initial condition $\u_0$. Since Problems
\ref{pb:maxET1D}, \ref{pb:maxET3D} and \ref{pb:PhiL4} are nonconvex,
their solutions discussed in Sections \ref{sec:1D} and \ref{sec:3D}
were obtained by {\em locally} maximizing $\varphi$ over flow
trajectories parameterized by the initial data $\u_0$ and hence may
not saturate the global maxima $\overline{\varphi}$.

On the other hand, it is possible to obtain upper bounds on the
supremum $\overline{\varphi}$ by exploiting the structure of the
governing equation, yet without reference to individual trajectories.
It has been shown in \cite{FantuzziGoluskin2020} that defining an
auxiliary function $V \; : \; \X \rightarrow \RR$ with the Lie
derivative {$\L V(\u(t)) := \big\langle \bnabla V(\u(t)),
  \f(\u(t)) \big\rangle_{\X} = dV(\u(t))/dt$,} such upper bounds can
be deduced by solving the following optimization problem
\begin{subequations} 
\label{eq:infV}
\begin{align}
\overline{\varphi} \ \le \ \inf_{V} \sup_{\u \in X} V(\u) &,  \label{eq:infVa}  \\
\L V(\u) & \le 0, \quad \u \in \X,  \label{eq:infVb} \\
\varphi(\u) - V(\u) & \le 0, \quad \u \in \X,  \label{eq:infVc} 
\end{align}
\end{subequations}
which is independent of any particular solution trajectories.
Importantly, in contrast to problem \eqref{eq:supphi}, {the
  outer minimization problem in \eqref{eq:infVa}} is convex and the
two problems are dual to each other \cite{l69} (under some additional
conditions this duality is strong).  In certain simple {cases
  problem \eqref{eq:infV} can be solved} analytically.  {Its
  numerical solution is in principle possible provided the inner
  maximization subproblem in \eqref{eq:infVa} can be suitably relaxed
  and the set of auxiliary functions $V$ is made finite-dimensional.
  For example,} when the auxiliary function $V(\u)$ and the function
$\f(\u)$ in the governing system are polynomial, inequality constraint
\eqref{eq:infVb} can be interpreted as imposing the non-negativity of
a polynomial which can then be expressed in terms of a sum of squares
(SoS) of some polynomial basis functions.  For PDE problems a
polynomial representation of $\f(\u)$ can be constructed using a
truncated Galerkin projection. These steps make it possible to
approximate problem \eqref{eq:infV} in terms of a semi-definite
optimization program for which many robust solution algorithms and
software packages are available.

The auxiliary function $V(\u)$ used in \eqref{eq:infV} is related to
the Lyapunov function employed in the study of nonlinear stability of
fixed points (except that, unlike the Lyapunov function, it need not
be positive semi-definite).  Formulations based on auxiliary functions
can also be used to obtain bounds on {infinite-time and space}
averages of various quantities of interest leading to convex
optimization problems analogous to \eqref{eq:infV}
\cite{hchgltf15,lhtch16,fghch16,Goluskin2018,TobascoGoluskinDoering2018}.
When the auxiliary function is fixed and quadratic whereas
optimization is performed with respect to the form of a certain
``background flow'', this bounding framework reduces to the background
method originally developed by Doering \& Constantin
\cite{DoeringConstantin1992} to obtain rigorous a priori bounds on
energy dissipation in wall-bounded flows. The background method has
been since used, both analytically and computationally, to derive
bounds on average quantities in different flows and we refer the
reader to \cite{Fantuzzi2021} for a recent survey of this topic.

In relation to the results reviewed in Sections \ref{sec:1D} and
\ref{sec:2D}, in \cite{FantuzziGoluskin2020} the authors used a
formulation based on auxiliary functions to rederive the a priori
bounds \eqref{eq:Etbound1D} and \eqref{eq:Pt2}, and to generalize the
former for the case of the ``fractional'' enstrophy $\E_{\alpha}
(u(t)) := \frac{1}{2}\int_0^1 \big| -\Delta^{\alpha / 2} u(t,x)\big|^2
\,dx$ relevant when the evolution is governed by the fractional
Burgers system \eqref{FBE}. In addition, by solving optimization
problem \eqref{eq:infV} for a Galerkin truncation of the Burgers
system \eqref{eq:Burgers} and a range of values of $\E_0$ they were
able to obtain upper bounds on $\E_T(\u_0)$ consistent with relation
\eqref{eq:maxETap}. This is a remarkable example of a situation when
the lower and upper bounds found by solving problems \eqref{eq:supphi}
and \eqref{eq:infV} coincide effectively closing the duality gap.
{Problems \ref{pb:maxET3D} and \ref{pb:PhiL4} can also be put in
  the framework of \eqref{eq:infV} and it is interesting to see
  whether it may be possible to develop suitable truncations and
  relaxations for the inner maximization problem what will yield
  computationally tractable semi-definite optimization programs.}

\section{Summary and Conclusions}
\label{sec:final}

In this paper we have presented a survey of recent progress in the
research program focused on a systematic computational search for
extreme behavior in different hydrodynamic models. Motivated by open
questions concerning the possibility of a finite-time blow-up in
solutions of the Navier-Stokes system \eqref{eq:NS} in 3D, these
investigations relied on solution of PDE-constrained optimization
problems with objective functionals chosen based on certain
conditional regularity results and a priori estimates available for
different models. {Families of local maximizers of} these PDE
optimization problems were {determined} numerically using
state-of-the-art adjoint-based gradient approaches formulated in the
continuous (infinite-dimensional) setting, cf.~Appendix \ref{sec:opt}.
We note that in addition to their numerous successful practical
applications involving, for example, shape optimization in
aerodynamics \cite{Jameson1988} and data assimilation in numerical
weather prediction \cite{k03}, similar optimization approaches have
also been employed in the study of some other fundamental problems in
fluid mechanics, namely, optimal mixing
\cite{EgglSchmid2020,MilesDoering2018}, transition to turbulence
\cite{rck12} and search for rare events in turbulence based on
instantons \cite{ggs15}.

\begin{table}
 
  \caption{Summary of a priori estimates considered in the research program discussed here 
together with information about their realizability.}
  \label{tab:estimates}
  \begin{center}
    \hspace*{-1.1cm}
    \begin{tabular}{lcc}      \hline
         \Bmp{2.5cm} \small \begin{center} {\sc Problem} \end{center} \Emp
      &  \Bmp{3.0cm} \small \begin{center} {\sc Estimate}  \end{center} \Emp   
      & \Bmp{3.5cm} \small \begin{center} {\sc Realizability }  \end{center} \Emp \\  
      \hline
      \Bmp{2.5cm}  \small {\begin{center} \smallskip 1D Burgers  \\ instantaneous \smallskip \end{center}} \Emp &  
      \small {$\frac{d\E}{dt} \leq \frac{3}{2}\left(\frac{1}{\pi^2\nu}\right)^{1/3}\E^{5/3}$}  & 
      \Bmp{3.5cm} \footnotesize {\begin{center} \smallskip {\sc Yes} \cite{ld08}  \smallskip  \end{center}} \Emp \\ 
      \hline 
      \Bmp{3.0cm} \small {\begin{center} \smallskip 1D Burgers  \\ finite-time \smallskip \end{center}} \Emp &  
      \small {$\max_{t \in [0,T]} \E(u(t)) \leq \left[\E_0^{1/3} + \frac{1}{16}\left(\frac{1}{\pi^2 \nu}\right)^{4/3}\E_0\right]^{3}$} &  
      \Bmp{3.5cm} \small {\begin{center} \smallskip {\sc No} \cite{ap11a,p12,p12b,FantuzziGoluskin2020} \smallskip  \end{center}} \Emp \\ 
      \hline 
      \Bmp{3.0cm} \small {\begin{center} \smallskip 2D Navier-Stokes  \\ instantaneous \smallskip\end{center}} \Emp &  
      \Bmp{7.0cm} \smallskip \centering \small
      \smallskip $\frac{d\P}{dt} \le C_2\sqrt{\log\left(\K^{1/2}/\nu\right)}\,  \P^{3/2}$ \smallskip \Emp& 
      \Bmp{3.5cm} \small {\begin{center} \smallskip {\sc Yes}  \cite{ap13a,ayala_doering_simon_2018}  \smallskip  \end{center}} \Emp  \\ 
      \hline 
      \Bmp{3.0cm} \small \begin{center} \smallskip 2D Navier-Stokes \\ finite-time \smallskip\end{center} \Emp &  
      \Bmp{7.0cm} \smallskip \centering \small 
      \smallskip $\max_{t \ge 0} \P(\u(t)) \le  \left(1 + \frac{a + b\sqrt{\ln Re_0 +c}}{4} Re_0 \right)^2 \P_0$ \smallskip \Emp  & 
      \Bmp{3.5cm} \small {\begin{center} \smallskip {\sc Yes}  \cite{ap13a,ayala_doering_simon_2018}  \smallskip  \end{center}} \Emp \\ 
      \hline 
      \Bmp{3.0cm} \small {\begin{center} \smallskip 3D Navier-Stokes  \\ instantaneous  \smallskip \end{center}} \Emp &  
      \small {$\frac{d\E}{dt} \leq \frac{27}{8\,\pi^4\,\nu^3} \E^3$} & \Bmp{3.5cm} \small {\begin{center} \smallskip {\sc Yes} \cite{ld08,ap16} \smallskip  \end{center}} \Emp  \\ 
      \hline 
      \Bmp{3.0cm} \small \begin{center} \smallskip\smallskip 3D Navier-Stokes  \\ finite-time \smallskip \end{center} \Emp &  
      \Bmp{7.0cm} \smallskip \centering \small  $\E(\u(t)) \le \frac{\E_0}{\sqrt{1 - 4 \frac{C \E_0^2}{\nu^3} t}}$ \\
      \smallskip $\int_0^T \| \u(\tau) \|_{L^4(\Omega)}^{8/3} \, d\tau \le  C \, \K_0^{4/3}$ \smallskip \Emp & \Bmp{3.5cm} \centering {\sc NO (???)} \cite{KangYunProtas2020,KangProtas2021} \smallskip \Emp \\
\hline
    \end{tabular}
\end{center}
\end{table}

The main results of the research program are summarized in Table
\ref{tab:estimates}. The main conclusion is that so far our search
based on the enstrophy and the Ladyzhenskaya-Prodi-Serrin conditions
\eqref{eq:supE} and \eqref{eq:LPS} has revealed no indication of
singularity formation in 3D Navier-Stokes flows with the optimal
initial conditions. However, the behavior exhibited by these different
extreme flows in terms of the maximum growth of enstrophy is in fact
similar, cf.~Figure \ref{fig:maxETvsE0}a, and analogous to what was
found in 1D Burgers flows, cf.~\eqref{eq:maxETap},
\eqref{eq:maxEt_vs_E0} and Figure \ref{fig:maxET1D}b. The results
discussed in Section \ref{sec:2D} demonstrated the sharpness of a
number of a priori estimates on the growth of palinstrophy in 2D, both
instantaneously and in finite time. We remark that the power-law
structure of estimates \eqref{eq:dEdt_estimate_E}, \eqref{eq:dPdt} and
of the empirical relation \eqref{eq:maxEt_vs_E0} can be justified with
simple arguments based on dimensional analysis
\cite{ld08,ayala_doering_simon_2018,KangYunProtas2020}. {The
  results surveyed here were obtained using different values of the
  viscosity coefficient $\nu$. To facilitate quantitative comparison
  between different problems, they can be rescaled to correspond to
  say $\nu = 1$ expressing the solution of \eqref{eq:NS} as $\u(t,\x)
  =: \nu \, \v(\nu t,\x)$, such that the rescaled solution $\v$ solves
  \eqref{eq:NS} with $\nu = 1$ and the time variable redefined as $\nu
  t$ (analogous approach applies to systems \eqref{eq:Burgers} and
  \eqref{eq:NS2D}).}

Somewhat paradoxically, the situation in 2D is more satisfactory than
in 1D where the key finite-time estimate \eqref{eq:Etbound1D} appears
not to be sharp, as indicated by the agreement of the results
discussed in Section \ref{sec:1D}, cf.~\eqref{eq:maxETap}, and the
upper bounds found in \cite{FantuzziGoluskin2020} by solving problem
\eqref{eq:infV}. Thus, rigorously improving this estimate remains an
open problem in PDE analysis and important progress towards this goal
has already been made in \cite{p12,p12b}. {More specifically, it
  is interesting to see whether the asymptotic estimate
  $\O(\E_0^{3/2})$ obtained in \cite{p12} for the maximum growth of
  enstrophy implies a rigorous upper bound on $\max_{t \ge 0} \E(t)$
  and whether the required assumptions on the regularity of the
  initial data can be weakened.}

Moving forward, the search for singular behavior in 3D Navier-Stokes
flows can be broadened by considering optimization problems analogous
to Problem \ref{pb:PhiL4} with objective functionals based on
conditional regularity results generalizing \eqref{eq:LPS} to include
norms of derivatives of different order of the velocity field
\cite{Gibbon2018}. In addition, this research program will be
broadened to include search for potential singularities in 3D Euler
flows which can be sought with similar approaches.

\appendix

\section{Solution of Optimization Problems}
\label{sec:opt}

In this appendix we provide some comments about the numerical
approaches employed to {find families of local maximizers in}
the optimization problems discussed in Sections \ref{sec:1D},
\ref{sec:2D} and \ref{sec:3D}. Since Problems
\ref{pb:maxdEdt1D}--\ref{pb:PhiL4} were designed to test certain
subtle properties of the underlying PDEs, we chose to formulate the
solution approaches in the continuous (``optimize-then-discretize'')
setting, where the optimality conditions, constraints and gradient
expressions are derived based on the original PDEs before being
discretized for the purpose of numerical evaluation, instead of the
alternative ``discretize-then-optimize'' approach often used in
applications \cite{g03}. In general, local maximizers in Problems
\ref{pb:maxdEdt1D}--\ref{pb:PhiL4} can be approximated using discrete
gradient flows with gradient expressions and constraints specific to
different problems. To fix attention, here we will describe in some
detail the approaches to solving Problems \ref{pb:maxdEdt3D} and
\ref{pb:maxET3D}, and then provide comment how to adapt them to solve
Problems \ref{pb:maxdPdt} and \ref{pb:PhiL4} (Problem
\ref{pb:maxdEdt1D} is solvable analytically, whereas Problem
\ref{pb:maxET1D} is a simpler 1D version of Problem \ref{pb:maxET3D}).
Finally, we will also provide some details about numerical
approximations.

\subsection{Solution of Problem \ref{pb:maxdEdt3D}}
\label{sec:optP4}

For a given value of $\E_0$, a local maximizer $\tuE$ of Problem
\ref{pb:maxdEdt3D} can be found as $\tuE = \lim_{n\rightarrow \infty}
\u_{\E_0}^{(n)}$ using the following iterative procedure representing
a discretization of a gradient flow projected on $\mathcal{S}_{\E_0}$
\begin{equation}
\begin{aligned}
\u_{\E_0}^{(n+1)} & =  \mathbb{P}_{\mathcal{S}_{\E_0}}\left(\;\u^{(n)}_{\E_0} + \tau_n \nabla\R\left(\u^{(n)}_{\E_0}\right)\;\right), \\ 
\u_{\E_0}^{(1)} & =  \u^0,
\end{aligned}
\label{eq:desc}
\end{equation}
where $\u^{(n)}_{\E_0}$ is an approximation of the maximizer obtained
at the $n$-th iteration, $\u^0$ is the initial guess and $\tau_n$ is
the length of the step in the direction of the gradient
$\nabla\R(\u^{(n)}_{\E_0})$.  Projection onto the constraint manifold
$\mathcal{S}_{\E_0}$ is performed using the composite operator
$\mathbb{P}_{\mathcal{S}_{\E_0}} \; : \; H^2(\Omega) \rightarrow
\mathcal{S}_{\E_0}$ defined as
\begin{subequations} 
\label{eq:Proj} 
\begin{align}
\mathbb{P}_{\mathcal{S}_{\E_0}}(\u) & = \P_{\E_0}\left(\Pi_0 (\u)\right), \qquad \text{where}  \label{eq:PP} \\
\Pi_0(\u) & = \u - \bnabla\left[\laplacian^{-1}(\bnabla\cdot\u)\right], 
\label{eq:Phodge} \\
\P_{\E_0}(\u) & = \sqrt{\frac{\E_0}{\E\left(\u\right)}}\,\u
\label{eq:Pnorm}
\end{align}
\end{subequations}
in which \eqref{eq:Phodge} and \eqref{eq:Pnorm} represent,
respectively, enforcement of the incompressibility condition and
normalization related to the enstrophy constraint.

A key step in procedure \eqref{eq:desc} is evaluation of the gradient
$\nabla\R(\u)$ of the objective functional $\R(\u)$, cf.
\eqref{eq:dEdt}, representing its (infinite-dimensional) sensitivity
to perturbations of the velocity field $\u$, and it is essential that
the gradient be characterized by the required regularity, namely,
$\nabla\R(\u) \in H^2(\Omega)$.  This is, in fact, guaranteed by the
Riesz representation theorem \cite{l69} applicable because the
G\^ateaux differential $\R'(\u;\cdot) : H^2(\Omega) \rightarrow \RR$,
defined as $\R'(\u;\u') := \lim_{\epsilon \rightarrow 0}
\epsilon^{-1}\left[\R(\u+\epsilon \u') - \R(\u)\right]$ for some
perturbation $\u' \in H^2(\Omega)$, is a bounded linear functional on
$H^2(\Omega)$.  The G\^ateaux differential can be computed directly to
give
\begin{equation}
\R'(\u;\u') = \int_{\Omega}\left[\u'\cdot\bnabla\u\cdot\laplacian\u + 
\u\cdot\bnabla\u'\cdot\laplacian\u + 
\u\cdot\bnabla\u\cdot\laplacian\u' \right]\,d\xvec 
-2\nu\int_{\Omega}\laplacian^2\u\cdot\u'\,d\xvec
\label{eq:dR}
\end{equation}
from which, by the Riesz representation theorem, we obtain
\begin{equation}
\R'(\u;\u') 
= \Big\langle \nabla\R(\u), \u' \Big\rangle_{H^2(\Omega)}
= \Big\langle \nabla^{L^2}\R(\u), \u' \Big\rangle_{L^2(\Omega)}
\label{eq:riesz}
\end{equation}
with the Riesz representers $\nabla\R(\u)$ and $\nabla^{L^2}\R(\u)$
being the gradients computed with respect to the $H^2$ and $L^2$
topology, respectively. We remark that, while the $H^2$ gradient is
used exclusively in the actual computations, cf.  \eqref{eq:desc}, the
$L^2$ gradient is computed first as an intermediate step.  Identifying
the G\^ateaux differential \eqref{eq:dR} with the $L^2$ inner product
and performing integration by parts yields
\begin{equation}
\nabla^{L^2}\R(\u) = \laplacian\left( \u\cdot\bnabla\u \right) + (\bnabla\u)^T\laplacian\u - 
\u\cdot\bnabla(\laplacian\u) - 2\nu\laplacian^2\u.
\label{eq:gradRL2}
\end{equation}
The inner product in $H^2(\Omega)$ is defined here as $\big\langle
\mathbf{z}_1, \mathbf{z}_2 \big\rangle_{H^2(\Omega)} := \int_{\Omega}
\mathbf{z}_1 \cdot \mathbf{z}_2 + \ell_1^2 \,\bnabla \mathbf{z}_1
\colon \bnabla \mathbf{z}_2 + \ell_2^4 \,\Delta \mathbf{z}_1 \cdot
\Delta \mathbf{z}_2 \, d\x$, $\forall\,\mathbf{z}_1, \mathbf{z}_2 \in
H^2(\Omega)$, where $\ell_1,\ell_2\in \RR_+$ are parameters with the
meaning of length scales (clearly, the inner products are equivalent
as long as $0 < \ell_1,\ell_2 < \infty$). Identifying the G\^ateaux
differential \eqref{eq:dR} with the $H^2$ inner product, integrating
by parts and using \eqref{eq:gradRL2}, we obtain the required $H^2$
gradient $\nabla\R$ as a solution of the elliptic boundary-value
problem
\begin{equation}
\begin{aligned}
&\left[ \Id \, - \,\ell_1^2 \,\Delta + \,\ell_2^4 \,\Delta^2 \right] \nabla\R
= \nabla^{L^2} \R  \qquad \text{in} \ \Omega, \\
& \text{Periodic Boundary Conditions}.
\end{aligned}
\label{eq:gradRH2}
\end{equation}
As shown in \cite{pbh04}, extraction of gradients in spaces of
smoother functions such as $H^2(\Omega)$ can be interpreted as
low-pass filtering of the $L^2$ gradients with parameters $\ell_1$ and
$\ell_2$ acting as the cut-off length-scales. The values of $\ell_1$
and $\ell_2$ can significantly affect the rate of convergence of the
iterative procedure \eqref{eq:desc}.

The step size $\tau_n$ in algorithm \eqref{eq:desc} is computed
as
\begin{equation}\label{eq:tau_n}
\tau_n = \mathop{\argmax}_{\tau>0} \left\{ \R\left[\mathbb{P}_{\mathcal{S}_{\E_0}}
\left( \;\u^{(n)} + \tau\,\nabla\R(\u^{(n)}) \;\right)\right] \right\}
\end{equation}
which is done using a suitable derivative-free line-search algorithm
\cite{r06}. Equation \eqref{eq:tau_n} can be interpreted as a
modification of a standard line search method where optimization is
performed following an arc (a geodesic) lying on the constraint
manifold $\mathcal{S}_{\E_0}$, rather than a straight line.

To ensure the maximizers $\tuE$ obtained for different values of
$\E_0$ lie on the same {maximizing} branch we use a continuation
approach, where the maximizer $\tuE$ is {employed} as the initial
guess $\u^0$ to compute $\widetilde{\mathbf{u}}_{\E_0+\Delta\E}$ using
\eqref{eq:desc} at the next enstrophy level for some sufficiently
small $\Delta\E > 0$.  We refer the reader to \cite{ap16} for further
details and add that in their seminal study \cite{ld08} Lu and Doering
used the alternative ``discretize-then-optimize'' approach.

In addition to some obvious simplifications, solution of Problem
\ref{pb:maxdPdt} does involve one important complication, namely, the
constraint manifold $\mathcal{W}_{\K_0,\P_0}$ is defined as an
intersection of two nonlinear manifolds. As a result, the projection
operator $\mathbb{P}_{\mathcal{W}_{\K_0,\P_0}}$ has a more complicated
structure: while the energy constraint
$\frac{1}{2}\int_\Omega|\bnabla\psi|^2\,d\Omega = \K_0$ is enforced
using normalization analogous to \eqref{eq:Pnorm}, the palinstrophy
constraint $\frac{1}{2}\int_\Omega|\bnabla\Delta\psi|^2\,d\Omega =
\P_0$ is satisfied by solving an inner optimization problem
$\min_{\phi \in H^4(\Omega)} \; (1/2)\left[ \P(\phi) - \P_0 \right]^2$
subject to $\K(\phi) = \K_0$ each time the objective functional is
evaluated in the discrete gradient flow \eqref{eq:desc}.

\subsection{Solution of Problem \ref{pb:maxET3D}}
\label{sec:optP5}

{Local maximizers of Problem \ref{pb:maxET3D} are determined}
with an approach similar to the method described in Section
\ref{sec:optP4} with one important difference, namely, the gradient
$\nabla\ET$ now needs to account for the flow evolution which is done
using methods of the adjoint calculus \cite{KangYunProtas2020}.  Given
the definition of the objective functional $\E_T(\u_0)$, its
G\^{a}teaux differential can be expressed as
\begin{equation}
\E'_T(\u_0;\u_0') = \int_\Omega (\rot\u(T,\x))\cdot(\rot\u'(T,\x))  \,d\x 
= \int_\Omega \bDelta\u(T,\x))\cdot\u'(T,\x)  \,d\x,
\label{eq:dET}
\end{equation}
where the last equality follows from integration by parts and the
vector identity $\bnabla\times(\bnabla\times\z) =
\bnabla(\bnabla\cdot\z) - \bDelta\z$, whereas the perturbation field
$\u' = \u'(t,\x)$ is a solution of the Navier-Stokes system linearized
around the trajectory corresponding to the initial data $\u_0$
\cite{g03}, i.e.,
\begin{subequations}
\label{eq:lNSE3D}
\begin{align}
 \L\begin{bmatrix} \u' \\ p' \end{bmatrix} := 
& \begin{bmatrix}
\partial_{t}\u'+\u'\cdot\bnabla\u+\u\cdot\bnabla\u'+\bnabla p'-\nu\bDelta\u' \\
\bnabla\cdot\u'
\end{bmatrix} = \begin{bmatrix} \mathbf{0} \\ 0\end{bmatrix}, \label{eq:lNSE3Da} \\
 \u'(0)= &\u_0' \label{eq:lNSE3Db}
\end{align}
\end{subequations}
which is subject to the periodic boundary conditions and where $p'$ is
the perturbation pressure.

We note that expression \eqref{eq:dET} for the G\^{a}teaux
differential is not consistent with the Riesz form \eqref{eq:riesz},
because the perturbation $\u_0'$ of the initial data does not appear
in it explicitly as a factor, but is instead hidden as the initial
{condition} in the linearized problem, cf.~\eqref{eq:lNSE3Db}. In
order to transform \eqref{eq:dET} to the Riesz form, we introduce the
{\em adjoint state} $\u^* \; : \; [0,T]\times\Omega \rightarrow
\RR^3$ and $p^* \; : \; [0,T]\times\Omega \rightarrow \RR$, and the
following duality-pairing relation
\begin{equation}
\begin{aligned}
\left( \L\begin{bmatrix} \u' \\ p' \end{bmatrix}, \begin{bmatrix} \u^* \\ p^* \end{bmatrix} \right)
:= & \int_0^T \int_{\Omega} \L\begin{bmatrix} \u' \\ p' \end{bmatrix} \cdot \begin{bmatrix} \u^* \\ p^* \end{bmatrix} \, d\x \, dt 
= \left( \begin{bmatrix} \u' \\ p' \end{bmatrix}, \L^*\begin{bmatrix} \u^* \\ p^* \end{bmatrix}\right) + \\
\phantom{=} & {\underbrace{\int_\Omega \u'(T,\x)\cdot\u^*(T,\x)  \,d\x}_{\E'_T(\u_0;\u_0')}} - 
\int_\Omega \u'(0,\x)\cdot\u^*(0,\x)  \,d\x = 0.
\end{aligned}
\label{eq:dual}
\end{equation}
Performing integration by parts with respect to {both space and time
  then} allows us to define the {\em adjoint system} as
\begin{subequations}
\label{eq:aNSE3D}
\begin{align}
 \L^*\begin{bmatrix} \u^* \\ p^* \end{bmatrix} := 
& \begin{bmatrix}
-\partial_{t}\u^*-\left[\bnabla\u^*+\left(\bnabla\u^{*}\right)^T\right]\u-\bnabla p^*-\nu\bDelta\u^* \\
-\bnabla\cdot\u^*
\end{bmatrix}  = \begin{bmatrix} {\mathbf{0}}  \\ 0\end{bmatrix}, \label{eq:aNSE3Da} \\
 \u^*(T)= & {\bDelta\u}  \label{eq:aNSE3Db}
\end{align}
\end{subequations}
which is also subject to the periodic boundary conditions. We note
that in identity \eqref{eq:dual} all boundary terms resulting from
integration by parts {with respect to the space variables} vanish due
to the periodic boundary conditions. The term $\int_\Omega
\u'(T,\x)\cdot\u^*(T,\x) \,d\x$ {resulting from integration by parts
  with respect to time is equal to the G\^ateaux differential
  \eqref{eq:dET} due to the judicious} choice of the terminal
condition \eqref{eq:aNSE3Db}, {such that} identity \eqref{eq:dual}
implies $\E'_T(\u_0;\u_0') = \int_\Omega \u'_0(\x)\cdot\u^*(0,\x)
\,d\x$, from which we deduce the following expression for the $L^2$
gradient
\begin{equation}
\nabla^{L^2}\ET = \u^*(0).
\label{eq:gradL2}
\end{equation}
The corresponding $H^1$ Sobolev gradient $\nabla\ET$ is then computed
as in Section \ref{sec:optP4}, using the Riesz identity to obtain an
elliptic boundary-value problem satisfied by the Sobolev gradient,
cf.~ \eqref{eq:riesz} and \eqref{eq:gradRH2}.

\subsection{Numerical Implementation}
\label{sec:numer}

Since Problems \ref{pb:maxdEdt1D}--\ref{pb:PhiL4} are all defined on
periodic domains, they can be accurately discretized in space using
standard Fourier pseudospectral methods with dealiasing
\cite{b01,canuto:SpecMthd}.  For the time-dependent problems, the time
discretization was performed using semi-implicit Runge-Kutta methods.
For 3D Problems \ref{pb:maxdEdt3D}, \ref{pb:maxET3D} and
\ref{pb:PhiL4} typical spatial resolutions varied from $128^3$ to
$512^3$ gridpoints which required massively parallel implementations
based on the Message Passing Interface (MPI). Solution of a single
instance of Problem \ref{pb:maxET3D} or \ref{pb:PhiL4} usually
required a computational time of $\O(10^2)$ hours on $\O(10^2)$ CPU
cores. The reader is referred to
\cite{ap11a,ap13a,ap16,KangYunProtas2020,KangProtas2021} for further
technical details.

\vskip6pt
\enlargethispage{20pt}


\dataccess{This article has no additional data.}



\funding{The author acknowledges the support through an NSERC (Canada)
  Discovery Grant.  Computational resources were provided by Compute
  Canada under its Resource Allocation Competitions.}

\ack{This work is dedicated to the memory of the late Charlie Doering,
  our dear friend and collaborator, who inspired us to pursue this
  research direction. The original plan was for this review article to
  be written jointly with Charlie, but his untimely death took him
  before writing could begin. The author thanks the current and former
  members of his research group: Diego Ayala, Di Kang, Pritpal ``Pip''
  Matharu, Diogo Po\c{c}as, Elkin Ramirez, Adam \'Sliwiak, Dongfang
  Yun and Xinyu Zhao for their contributions to the research program
  surveyed in this paper. The author also thanks Miguel Bustamante,
  Sergei Chernyshenko, Giovanni Fantuzzi, David Goluskin, John Gibbon,
  Thomas Y.~Hou, Evan Miller, Koji Ohkitani, Dmitry Pelinovsky and
  Tsuyoshi Yoneda for many enlightening and enjoyable discussions.
  {The author acknowledges useful feedback provided by the
    referees.}}



\vskip2pc













\begin{thebibliography}{99}

\bibitem{d09}
Doering CR. 2009  The {3D Navier-Stokes} Problem. {\em Annual Review of Fluid
  Mechanics} \textbf{41}, 109--128.

\bibitem{Robinson2020}
Robinson JC. 2020  {The Navier-Stokes regularity problem}. {\em {Phil. Trans.
  R. Soc. A}} \textbf{378}, 20190526.

\bibitem{kl04}
Kreiss H, Lorenz J. 2004 {\em Initial-Boundary Value Problems and the
  {Navier-Stokes} Equations} vol.~47{\em Classics in Applied Mathematics}.
SIAM.

\bibitem{f00}
Fefferman CL. 2000  Existence and Smoothness of the {Navier-Stokes} Equation.
  available at
  {\tt{http://www.claymath.org/sites/default/files/navierstokes.pdf}}.
{Clay Millennium Prize Problem Description}.

\bibitem{l34}
Leray J. 1934  Sur le mouvement d'un liquide visqueux emplissant l'espace. {\em
  Acta Mathematica} \textbf{63}, 193--248.

\bibitem{BuckmasterVicol2019}
Buckmaster T, Vicol V. 2019  {Nonuniqueness of weak solutions to the
  Navier-Stokes equation}. {\em Annals of Mathematics} \textbf{189}, 101--144.

\bibitem{gbk08}
Gibbon JD, Bustamante M, Kerr RM. 2008  The three--dimensional {Euler}
  equations: singular or non--singular?. {\em Nonlinearity} \textbf{21},
  123--129.

\bibitem{ap11a}
Ayala D, Protas B. 2011  On Maximum Enstrophy Growth in a Hydrodynamic System.
  {\em Physica D} \textbf{240}, 1553--1563.

\bibitem{ap13a}
Ayala D, Protas B. 2014a  Maximum Palinstrophy Growth in {2D} Incompressible
  Flows. {\em Journal of Fluid Mechanics} \textbf{742}, 340--367.

\bibitem{ap13b}
Ayala D, Protas B. 2014b  Vortices, maximum growth and the problem of
  finite-time singularity formation. {\em Fluid Dynamics Research} \textbf{46},
  031404.

\bibitem{ap16}
Ayala D, Protas B. 2017  Extreme Vortex States and the Growth of Enstrophy in
  {3D} Incompressible Flows. {\em Journal of Fluid Mechanics} \textbf{818},
  772--806.

\bibitem{Yun2018}
Yun D, Protas B. 2018  {Maximum Rate of Growth of Enstrophy in Solutions of the
  Fractional Burgers Equation}. {\em Journal of Nonlinear Science} \textbf{28},
  395--422.

\bibitem{KangYunProtas2020}
Kang D, Yun D, Protas B. 2020  {Maximum amplification of enstrophy in
  three-dimensional Navier-Stokes flows}. {\em Journal of Fluid Mechanics}
  \textbf{893}, A22.

\bibitem{af05}
Adams RA, Fournier JF. 2005 {\em {Sobolev} Spaces}.
Elsevier.

\bibitem{RobinsonRodrigoSadowski2016}
Robinson JC, Rodrigo JL, Sadowski W. 2016 {\em {The Three-Dimensional
  Navier-Stokes Equations: Classical Theory}}.
Cambridge University Press.

\bibitem{ft89}
Foias C, Temam R. 1989  Gevrey Class Regularity for the Solutions of the
  {Navier--Stokes} Equations. {\em Journal of Functional Analysis} \textbf{87},
  359--369.

\bibitem{KisLad57}
Kiselev AA, Ladyzhenskaya OA. 1957  On the existence and uniqueness of the
  solution of the nonstationary problem for a viscous, incompressible fluid.
  {\em Izv. Akad. Nauk SSSR Ser. Mat} \textbf{21}, 655--680.

\bibitem{Prodi1959}
Prodi G. 1959  {{Un teorema di unicit{\`a} per le equazioni di Navier-Stokes}}.
  {\em Annali di Matematica Pura ed Applicata} \textbf{48}, 173--182.

\bibitem{Serrin1962}
Serrin J. 1962  {On the interior regularity of weak solutions of the
  Navier-Stokes equations}. {\em Archive for Rational Mechanics and Analysis}
  \textbf{9}, 187--195.

\bibitem{Gibbon2018}
Gibbon JD. 2018  {Weak and Strong Solutions of the 3D Navier--Stokes Equations
  and Their Relation to a Chessboard of Convergent Inverse Length Scales}. {\em
  Journal of Nonlinear Science}.
(published on-line).

\bibitem{Escauriaza2003}
Escauriaza L, Seregin GA, Sverak V. 2003  L3,$\infty$-solutions of the
  Navier-Stokes equations and backward uniqueness. {\em Russian Mathematical
  Surveys} \textbf{58}, 211--250.

\bibitem{Tao2020}
Tao T. 2020  {Quantitative bounds for critically bounded solutions to the
  Navier-Stokes equations}. arXiv:1908.04958.

\bibitem{Constantin1991}
Constantin P. 1991  {Remarks on the Navier-Stokes equations}. In Sirovich L,
  editor, {\em New Perspectives in Turbulence} pp. 229--261. Berlin: Springer.

\bibitem{KangProtas2021}
Kang D, Protas B. 2021  {Searching for Singularities in Navier-Stokes Flows
  Based on the Ladyzhenskaya-Prodi-Serrin Conditions}. arXiv:2110.06130.

\bibitem{bkm84}
Beale JT, Kato T, Majda A. 1984  {Remarks on the breakdown of smooth solutions
  for the $3$-D Euler equations}. {\em Comm. Math. Phys.} \textbf{94}, 61--66.

\bibitem{Constantin1986}
Constantin P. 1986  Note on loss of regularity for solutions of the 3---D
  incompressible euler and related equations. {\em Communications in
  Mathematical Physics} \textbf{104}, 311--326.

\bibitem{ElgindiJeong2018}
Elgindi TM, Jeong IJ. 2019  {Finite-Time Singularity Formation for Strong
  Solutions to the Axi-symmetric 3D Euler Equations}. {\em Annals of PDE}
  \textbf{5}, 16.

\bibitem{ld08}
Lu L, Doering CR. 2008  Limits on Enstrophy Growth for Solutions of the
  Three-dimensional {Navier--Stokes} Equations. {\em Indiana University
  Mathematics Journal} \textbf{57}, 2693--2727.

\bibitem{Giga1986}
Giga Y. 1986  {Solutions for semilinear parabolic equations in $L^p$ and
  regularity of weak solutions of the Navier-Stokes system}. {\em Journal of
  Differential Equations} \textbf{62}, 186--212.

\bibitem{RobinsonSadowskiSilva2012}
Robinson JC, Sadowski W, Silva RP. 2012  Lower bounds on blow up solutions of
  the three-dimensional Navier~Stokes equations in homogeneous Sobolev spaces.
  {\em Journal of Mathematical Physics} \textbf{53}, 115618.

\bibitem{RobinsonSadowski2014}
Robinson JC, Sadowski W. 2014  A local smoothness criterion for solutions of
  the 3D Navier-Stokes equations. {\em Rendiconti del Seminario Matematico
  della Universit\`a di Padova} \textbf{131}, 159--178.

\bibitem{bmonmu83}
Brachet ME, Meiron DI, Orszag SA, Nickel BG, Morf RH, Frisch U. 1983
  Small-scale Structure of the {Taylor-Green} Vortex. {\em Journal of Fluid
  Mechanics} \textbf{130}, 411--452.

\bibitem{ps90}
Pumir A, Siggia E. 1990  Collapsing solutions to the {3D Euler} equations. {\em
  Phys. Fluids A} \textbf{2}, 220--241.

\bibitem{b91}
Brachet ME. 1991  Direct Simulation of Three-dimensional Turbulence in the
  {Taylor-Green} Vortex. {\em Fluid Dynamics Research} \textbf{8}, 1--8.

\bibitem{k93}
Kerr RM. 1993  Evidence for a Singularity of the Three-dimensional,
  Incompressible {Euler} Equations. {\em Phys. Fluids A} \textbf{5},
  1725--1746.

\bibitem{p01}
Pelz RB. 2001  Symmetry and the hydrodynamic blow-up problem. {\em Journal of
  Fluid Mechanics} \textbf{444}, 299--320.

\bibitem{bk08}
Bustamante MD, Kerr RM. 2008  {3D Euler about a 2D symmetry plane}. {\em
  Physica D} \textbf{237}, 1912--1920.

\bibitem{oc08}
Ohkitani K, Constantin P. 2008  {Numerical study of the Eulerian--Lagrangian
  analysis of the {Navier-Stokes} turbulence}. {\em Phys. Fluids} \textbf{20},
  1--11.

\bibitem{o08}
Ohkitani K. 2008  A miscellany of basic issues on incompressible fluid
  equations. {\em Nonlinearity} \textbf{21}, 255--271.

\bibitem{ghdg08}
Grafke T, Homann H, Dreher J, Grauer R. 2008  Numerical simulations of possible
  finite-time singularities in the incompressible {Euler} equations: comparison
  of numerical methods. {\em Physica D} \textbf{237}, 1932--1936.

\bibitem{h09}
Hou TY. 2009  Blow-up or No Blow-up? A Unified Computational and Analytic
  Approach to {3D} Incompressible {Euler} and {Navier--Stokes} Equations. {\em
  Acta Numerica} pp. 277--346.

\bibitem{opc12}
Orlandi P, Pirozzoli S, Carnevale GF. 2012  Vortex events in {Euler} and
  {Navier-Stokes} simulations with smooth initial conditions. {\em Journal of
  Fluid Mechanics} \textbf{690}, 288--320.

\bibitem{bb12}
Bustamante MD, Brachet M. 2012  Interplay between the {Beale-Kato-Majda}
  theorem and the analyticity-strip method to investigate numerically the
  incompressible {Euler} singularity problem. {\em Phys. Rev. E} \textbf{86},
  066302.

\bibitem{opmc14}
Orlandi P, Pirozzoli S, Bernardini M, Carnevale GF. 2014  A minimal flow unit
  for the study of turbulence with passive scalars. {\em Journal of Turbulence}
  \textbf{15}, 731--751.

\bibitem{CampolinaMailybaev2018}
Campolina CS, Mailybaev AA. 2018  Chaotic Blowup in the 3D Incompressible Euler
  Equations on a Logarithmic Lattice. {\em Phys. Rev. Lett.} \textbf{121},
  064501.

\bibitem{dggkpv13}
Donzis DA, Gibbon JD, Gupta A, Kerr RM, Pandit R, Vincenzi D. 2013  Vorticity
  moments in four numerical simulations of the {3D Navier-Stokes} equations.
  {\em Journal of Fluid Mechanics} \textbf{732}, 316--331.

\bibitem{k13}
Kerr RM. 2013  Swirling, turbulent vortex rings formed from a chain reaction of
  reconnection events. {\em Physics of Fluids} \textbf{25}, 065101.

\bibitem{gdgkpv14}
Gibbon J, Donzis D, Gupta A, Kerr R, Pandit R, Vincenzi D. 2014  {Regimes of
  nonlinear depletion and regularity in the 3D Navier-Stokes equations}. {\em
  Nonlinearity} \textbf{27}.

\bibitem{k13b}
Kerr RM. 2013  Bounds for {Euler} from vorticity moments and line divergence.
  {\em Journal of Fluid Mechanics} \textbf{729}, R2.

\bibitem{Kerr2018}
Kerr RM. 2018  {Enstrophy and circulation scaling for Navier-Stokes
  reconnection}. {\em Journal of Fluid Mechanics} \textbf{839}, R2.

\bibitem{MoffattKimura2019a}
Moffatt HK, Kimura Y. 2019a  {Towards a finite-time singularity of the
  Navier-Stokes equations Part 1. Derivation and analysis of dynamical system}.
  {\em Journal of Fluid Mechanics} \textbf{861}, 930--967.

\bibitem{MoffattKimura2019b}
Moffatt HK, Kimura Y. 2019b  Towards a finite-time singularity of the
  {Navier-Stokes} equations. Part 2. Vortex reconnection and singularity
  evasion. {\em Journal of Fluid Mechanics} \textbf{870}, R1.

\bibitem{mbf08}
Matsumoto T, Bec J, Frisch U. 2008  Complex-space singularities of {2D Euler}
  flow in Lagrangian coordinates. {\em Physica D} \textbf{237}, 1951--1955.

\bibitem{sc09}
Siegel M, Caflisch RE. 2009  Calculation of complex singular solutions to the
  {3D} incompressible {Euler} equations. {\em Physica D} \textbf{238},
  2368--2379.

\bibitem{lh14a}
Luo G, Hou TY. 2014a  {Potentially Singular Solutions of the 3D Axisymmetric
  Euler Equations}. {\em Proceedings of the National Academy of Sciences}
  \textbf{111}, 12968--12973.

\bibitem{lh14b}
Luo G, Hou TY. 2014b  {Toward the Finite-Time Blowup of the 3D Incompressible
  Euler Equations: a Numerical Investigation}. {\em SIAM: Multiscale Modeling
  and Simulation} \textbf{12}, 1722--1776.

\bibitem{HouHuang2021}
Hou TY, Huang D. 2021  {Potential Singularity Formation of 3D Axisymmetric
  Navier-Stokes Equations with Degenerate Diffusion Coefficients}.
  arXiv:2102.06663.

\bibitem{control:abergel1}
Abergel F, Temam R. 1990  On Some Control Problems in Fluid Mechanics. {\em
  Theoretical and Computational Fluid Dynamics} \textbf{1}, 303--325.

\bibitem{Fursikov2000}
Fursikov AV. 2000 {\em Optimal Control of Distributed Systems. Theory and
  Applications}.
Translations of Mathematical Monographs. American Mathematical Society.

\bibitem{g03}
Gunzburger MD. 2003 {\em Perspectives in Flow Control and Optimization}.
SIAM.

\bibitem{Troltzsch2010}
Tr{\"{o}}ltzsch F. 2010 {\em Optimal Control of Partial Differential Equations:
  Theory, Methods and Applications} vol. 112{\em Graduate Studies in
  Mathematics}.
American Mathematical Society.

\bibitem{bk07}
Bec J, Khanin K. 2007  Burgers turbulence. {\em Physics Reports} \textbf{447},
  1--66.

\bibitem{Biryuk2001}
Biryuk A{\`E}. 2001  Spectral Properties of Solutions of the Burgers Equation
  with Small Dissipation. {\em Functional Analysis and Its Applications}
  \textbf{35}, 1--12.

\bibitem{p12}
Pelinovsky D. 2012a  Sharp bounds on enstrophy growth in the viscous {Burgers}
  equation. {\em Proceedings of Royal Society A} \textbf{468}, 3636--3648.

\bibitem{p12b}
Pelinovsky D. 2012b  Enstrophy growth in the viscous {Burgers} equation. {\em
  Dynamics of Partial Differential Equations} \textbf{9}, 305--340.

\bibitem{FantuzziGoluskin2020}
Fantuzzi G, Goluskin D. 2020  Bounding Extreme Events in Nonlinear Dynamics
  Using Convex Optimization. {\em SIAM Journal on Applied Dynamical Systems}
  \textbf{19}, 1823--1864.

\bibitem{f15}
Flandoli F. 2015 {\em {Random Perturbation of PDEs and Fluid Dynamic Models}}.
Lecture Notes in Mathematics. Springer.

\bibitem{PocasProtas2018}
Po\c{c}as D, Protas B. 2018  {Transient growth in stochastic Burgers flows}.
  {\em Discrete \& Continuous Dynamical Systems --- B} \textbf{23}, 2371.

\bibitem{kns08}
Kiselev A, Nazaraov F, Shterenberg R. 2008  {Blow up and regularity for fractal
  Burgers equation}. {\em Dynamics of Partial Differential Equations}
  \textbf{5}, 211--240.

\bibitem{kp12}
Katz N, Pavlovi{\'{c}} N. 2002  A cheap {Caffarelli-Kohn-Nirenberg inequality
  for the Navier-Stokes} equation with hyper-dissipation. {\em Geometric {\&}
  Functional Analysis GAFA} \textbf{12}, 355--379.

\bibitem{Boritchev2018}
Boritchev A. 2018  {Decaying turbulence for the fractional subcritical Burgers
  equation}. {\em Discrete \& Continuous Dynamical Systems} \textbf{38},
  2229--2249.

\bibitem{RamirezProtas2021}
Ram\'irez E, Protas B. 2021  {Singularity Formation in the Deterministic and
  Stochastic Fractional Burgers Equation}. arXiv:2104.10759.

\bibitem{ayala_doering_simon_2018}
Ayala D, Doering CR, Simon TM. 2018  {Maximum palinstrophy amplification in the
  two-dimensional Navier-Stokes equations}. {\em Journal of Fluid Mechanics}
  \textbf{837}, 839--857.

\bibitem{td06}
Tran CV, Dritschel DG. 2006  Vanishing enstrophy dissipation in two-dimensional
  {Navier--Stokes} turbulence in the inviscid limit. {\em Journal of Fluid
  Mechanics} \textbf{559}, 107--116.

\bibitem{dfj10}
Dascaliuc R, Foias C, Jolly MS. 2010  Estimates on enstrophy, palinstrophy, and
  invariant measures for {2D} turbulence. {\em Journal of Differential
  Equations} \textbf{248}, 792--819.

\bibitem{JeongYoneda2021}
Jeong IJ, Yoneda T. 2021  {Enstrophy dissipation and vortex thinning for the
  incompressible 2D Navier{\textendash}Stokes equations}. {\em Nonlinearity}
  \textbf{34}, 1837--1853.

\bibitem{Sliwiak2017}
{\'S}liwiak A. 2017  {Maximum Rate of Growth of Enstrophy in the Navier-Stokes
  System on 2D Bounded Domains}. Master's thesis McMaster University.

\bibitem{tg37}
Taylor GI, Green AE. 1937  Mechanism of the production of small eddies from
  large ones. {\em Proceedings of the Royal Society of London A} \textbf{158},
  499--521.

\bibitem{pbh04}
Protas B, Bewley T, Hagen G. 2004  A comprehensive framework for the
  regularization of adjoint analysis in multiscale {PDE} systems. {\em Journal
  of Computational Physics} \textbf{195}, 49--89.

\bibitem{protas2008}
Protas B. 2008  Adjoint-Based optimization of {PDE} systems with alternative
  gradients. {\em Journal of Computational Physics} \textbf{227}, 6490--6510.

\bibitem{l69}
Luenberger D. 1969 {\em Optimization by Vector Space Methods}.
John Wiley and Sons.

\bibitem{hchgltf15}
Huang D, Chernyshenko S, Goulart P, Lasagna D, Tutty O, Fuentes F. 2015
  Sum-of-squares of polynomials approach to nonlinear stability of fluid flows:
  an example of application. {\em Proceedings of the Royal Society of London.
  Series A, Mathematical and physical sciences} \textbf{471}, 1--18.

\bibitem{lhtch16}
Lasagna D, Huang D, Tutty O, Chernyshenko S. 2016  Sum-of-Squares approach to
  feedback control of laminar wake flows. {\em Journal of Fluid Mechanics}
  \textbf{809}, 628--663.

\bibitem{fghch16}
Fantuzzi G, Goluskin D, Huang D, Chernyshenko SI. 2016  Bounds for
  Deterministic and Stochastic Dynamical Systems using Sum-of-Squares
  Optimization. {\em SIAM Journal on Applied Dynamical Systems} \textbf{15},
  1962--1988.

\bibitem{Goluskin2018}
Goluskin D. 2018  {Bounding Averages Rigorously Using Semidefinite Programming:
  Mean Moments of the Lorenz System}. {\em Journal of Nonlinear Science}
  \textbf{28}, 621--651.

\bibitem{TobascoGoluskinDoering2018}
Tobasco I, Goluskin D, Doering CR. 2018  Optimal bounds and extremal
  trajectories for time averages in nonlinear dynamical systems. {\em Physics
  Letters A} \textbf{382}, 382--386.

\bibitem{DoeringConstantin1992}
Doering CR, Constantin P. 1992  Energy dissipation in shear driven turbulence.
  {\em Phys. Rev. Lett.} \textbf{69}, 1648--1651.

\bibitem{Fantuzzi2021}
Fantuzzi G, Arslan A, Wynn A. 2021  {The background method: Theory and
  computations}. arXiv:2107.11206.

\bibitem{Jameson1988}
Jameson A. 1988  Aerodynamic design via control theory. {\em Journal of
  Scientific Computing} \textbf{3}, 233--260.

\bibitem{k03}
Kim N. 2003  {Remarks for the axisymmetric Navier-Stokes equations}. {\em
  Journal of Differential Equations} \textbf{187}, 226--239.

\bibitem{EgglSchmid2020}
Eggl MF, Schmid PJ. 2020  Mixing enhancement in binary fluids using optimised
  stirring strategies. {\em Journal of Fluid Mechanics} \textbf{899}, A24.

\bibitem{MilesDoering2018}
Miles CJ, Doering CR. 2018  A Shell Model for Optimal Mixing. {\em Journal of
  Nonlinear Science} \textbf{28}, 2153--2186.

\bibitem{rck12}
Rabin SME, Caulfield CP, Kerswell RR. 2012  Variational identification of
  minimal seeds to trigger transition in plane {Couette} flow. {\em Journal of
  Fluid Mechanics} \textbf{712}, 244--272.

\bibitem{ggs15}
Grafke T, Grauer R, Sch{\"a}fer T. 2015  The instanton method and its numerical
  implementation in fluid mechanics. {\em Journal of Physics A: Mathematical
  and Theoretical} \textbf{48}, 333001.

\bibitem{r06}
Ruszczy\'nski A. 2006 {\em Nonlinear Optimization}.
Princeton University Press.

\bibitem{b01}
Boyd JP. 2001 {\em Chebyshev and Fourier Spectral Methods}.
Dover.

\bibitem{canuto:SpecMthd}
Canuto C, Quarteroni A, Hussaini Y, Zang TA. 2006 {\em Spectral Methods}.
Scientific Computation. Springer.

\end{thebibliography}

\end{document}